\documentclass[a4paper,9pt]{article}

\usepackage[utf8]{inputenc}
\usepackage[english]{babel} 
\usepackage[T1]{fontenc}
\usepackage{amsmath}
\usepackage{amsfonts}
\usepackage{graphicx}
\usepackage{amssymb}
\usepackage{rotating}
\usepackage{color}
\usepackage[singlespacing]{setspace}
\usepackage[left=2.5cm, right=2.5cm, top=2.5cm, bottom=2.5cm]{geometry}

\newtheorem{thm}{Theorem}[section] 
\newtheorem{prop}[thm]{Proposition}
\newtheorem{cor}[thm]{Corollary}
\newtheorem{theo}[thm]{Theorem}
\newtheorem{defi}[thm]{Definition}
\newtheorem{lem}[thm]{Lemma}

\begin{document}
\newgeometry{left=2.2cm, right=2.2cm, top=2.8cm, bottom=2.8cm}
\begin{center}
\large{\textbf{GENERALIZATION OF ARNOLD'S $J^+$-INVARIANT FOR PAIRS OF IMMERSIONS}}\\
\textnormal{Hanna Häußler, 28.08.2024} \\
\textnormal{hanna.haeussler@uni-a.de}
\end{center}
\ \\
ABSTRACT: This paper introduces the $J^{2+}$-invariant for oriented pairs of generic immersions. This invariant behaves like Arnold's $J^+$-invariant for generic immersions as it is invariant when going through inverse tangencies and triple points, but changes when traversing direct tangencies. It has several useful properties, for example its independence of the $J^{+}$-invariants of the single immersions forming the pair. Also it is invariant under simultaneous orientation change. Therefore, one can define two $J^{2+}$-invariants for each pair depending on its orientation, those two invariants are not independent from each other. Furthermore the invariant is extended to the $J^{n+}$-invariant for links of $n$ oriented immersions.\\
\ \\
KEYWORDS: Arnold's $J^+$-invariant, Viro's formula, pairs of generic immersions, regularization of collisions\\
\pagenumbering{roman}
\tableofcontents 
\pagenumbering{arabic}
\ \\
\begin{center}\textbf{Introduction}\end{center}
\ \\
According to the Whitney-Graustein theorem [7], the rotation number of an immersion of the circle in the plane fully classifies such loops up to regular homotopy. During a generic regular homotopy, three critical events can occur: triple intersections, direct self-tangencies, and inverse self-tangencies. In [1] Arnold introduced three invariants $J^+,$ $J^-$ and $St$, which remain constant under two of the three events and change with one critical scenario. The focus of this paper is on the $J^+$-invariant, which remains unaffected by triple intersections and inverse self-tangencies but is sensitive to direct self-tangencies. \\
In many physical systems, for example those governed by conservative forces, self-intersections in orbits are not physically realizable. For example, the trajectory of a particle under a conservative force field in a plane typically cannot cross itself due to conservation of energy. Whereas, inverse self-tangencies can occur if the force depends on the velocity, as in the case of the Coriolis or Lorentz forces.\\
To analyze the movement of a particle in a conservative force field on the plane, Cieliebak-Frauenfelder-van Koert [2] introduce Stark-Zeeman systems, which model electron motion in electric and magnetic fields, and also describe the restricted three-body problem in celestial mechanics. In Stark-Zeeman systems, homotopies of periodic orbits can go through two additional critical events: (1) the creation or destruction of an external loop when the orbit's velocity temporarily reaches zero, and (2) the creation or destruction of a loop around the origin during a collision. When regularizing collisions at the origin pairs of immersions can occur as the preimage of an immersions, for example under the Levi-Civita, but also under the Birkhoff-regularization. \\
In this paper we introduce a $J^+$-like invariant for pairs and later also for $n$-links of immersions. Pairs of immersions can occur as preimages under regularization of collisions, but also if we consider the setting that we want a satellite to change from one orbit to another one. At best, these orbits intersect each other so that with a boost into the right direction at the intersection point, the satellite changes its orbit. \\
The $J^+$-like invariant for links of orbits can be used to study the problem whether different links can be connected by a family of non-bifurcating links. The $J^+$-invariant of the single immersions forming the link is not sufficient to study this question, as it doesn't see critical events between different periodic orbits. Thus periodic orbits in links can have the same $J^+$-invariants, but can still not be connected by a family of links of periodic orbits if the $J^{2+}$-invariants are different. Thus, the invariant helps us to get a better understanding of on the global picture how families of orbits are linked to each other. \\
Arnold defined his invariant by its behavior under homotopy and by associating a specific value of the invariant to the so called standard immersions $K_i$ with winding number $i$. Using Viro's formula [6], the invariant of arbitrary immersions in the plane is easy to compute. It is given through
\begin{equation}
J^+(S) = 1 + n - \sum_{C\subset \mathbb{C}\backslash S} \omega_C(S)^2 + \sum_p ind_p(S)^2 ,\nonumber
\end{equation} 
where $n$ is the amount of double points of the generic immersion $S$ and $\omega_C(S)$ is the winding number of the component $C$ of $\mathbb{C}\setminus \{S\}$. The index $ind_p(S)$ for all double points $p$ is defined as the arithmetic mean of the winding numbers of the four adjacent components. \\
\ \\
Before introducing a generalization of $J^+$, the paper gives an alternative proof for Viro's formula. This proof is straight forward and also similar to the approach in the standard proof by  Viro, but does not use $Rohklin's$ $type$ $formula$. It is included here, as it gives an intuition in how to adapt the formula to pairs of immersions.\\
Then, we introduce the $J^{2+}$-invariant for oriented pairs of generic immersions $K$ out of immersions $S_1$ and $S_2$ (we write $K(S_1,S_2)$)
\begin{equation}
J^{2+}(K) := 2 + n - \sum_{C \subset \mathbb{C}\setminus K} \omega_c(K)^2 + \sum_p ind_p(K)^2 + u(K)^2, \nonumber
\end{equation} 
where $n$, $\omega_c(K)$ and $ind_p(K)$ are defined analogously to Viro's formula for a pair of immersions $K$. The additional variable $u(K)$, called the encircling index, describes the relation between $S_1$ and $S_2$. \\
\ \\
\textbf{Theorem A}: \textit{The $J^{2+}$-invariant is constant under inverse tangency and triple points and increases (decreases) by the value of $2$ when going positively (negatively) through direct tangencies..} \\
\ \\
Then following theorem discusses the independency of $J^{2+}(K)$ and the $J^+$-invariants of $S_1$ and $S_2$.\\
\ \\
\textbf{Theorem B}: \textit{For arbitrary even numbers $x$, $y$ and $z$, there always exists an oriented pair of generic immersions $K(S_1,S_2)$ such that ${J}^{+}(S_1)=x$, ${J}^{+}(S_2)=y$ and ${J}^{2+}(K)=z$.} \\
\ \\
The proof is given through an algorithm which constructs suitable pairs of immersions satisfying the given conditions.\\
\ \\
As the $J^+$-invariant is independent of orientation, whereas the $J^{2+}$-invariant is not, but remains constant under simultaneous change of orientation, we introduce $\overline{K}$ as $K$ with changed orientation of one immersion. We study the relation of their $J^{2+}$-invariants and prove  that they are not independent of each other. \\
\ \\
\textbf{Theorem C:} \textit{Let $K$ be an arbitrary oriented pair out of the immersions $S_1$ and $S_2$, then the invariants  $J^+(S_1)=x$, $J^+(S_2)=y$, $J^{2+}(K)=z$ and $J^{2+}(\overline{K})=\overline{z}$ are even numbers satisfying $z + \overline{z} \geq 2 x + 2 y$.}\\
\ \\
Furthermore we give algorithms to construct pairs of immersions for all possible combinations of $J^+$- and $J^{2+}$-invariants.\\
\ \\
\textbf{Theorem D}: \textit{For arbitrary even numbers $x$, $y$ and $z$ and $\overline{z}$ satisfying $z + \overline{z} \geq 2 x + 2 y$, there always exists an oriented pair of generic immersions $K(S_1,S_2)$  such that ${J}^{+}(S_1)=x$, ${J}^{+}(S_2)=y$ and ${J}^{2+}(K)=z$ and ${J}^{2+}(\overline{K})=\overline{z}$.} \\
\ \\
We extend the invariant tolinks of $n$ immersions $(S_1,S_2,..,S_n)$ $\forall n \geq 3$ by
\begin{equation*}
    J^{n+} (S_1, S_2,..,S_n) = \sum_{i<j} \Bigl( J^{2+}(S_i, S_j)-J^+(S_i) - J^+(S_j)\Bigr) + \sum_{i=1}^n J^+(S_i) .
\end{equation*} 
\textbf{Theorem E}: \textit{The  $J^{n+}$-invariant for links of $n$ oriented immersions $S_1,S_2,..,S_n$ changes by the value of $2$ when going through direct tangency and remains constant under inverse tangency and triple points.}\\
\ \\
Also we can define an extension of the $J^-$-invariant which is constant under direct tangencies and triple points and decreases by the value of $2$ when going through positive inverse tangency. \\
\ \\
This paper has the following structure: In the first chapter, we recall the definition of Arnold's $J^+$-invariant and Viro's formula. An alternative proof of the formula is given in Section 1.2. The second chapter introduces the $J^{2+}$-invariant for oriented pairs of immersions and Theorem A is proven. In Section 2.2, we look at additional properties of the $J^{2+}$-invariant and prove Theorem B. In chapter 3, $\overline{K}$ and its $J^{2+}$-invariant is introduced and its dependency on $J^{2+}(K)$ and the $J^+$-invariants forming the pair are shown by proving Theorem C. Theorem D tells us that we can construct suitable pairs of immersions for given arbitrary values fulfilling the necessary constraints. Section 3.3 introduces the $J^{2-}$-invariant and studies its relation to $J^{2+}(\overline{K})$. In Chapter 4 we extend the invariant to a system of $n$ immersions with the introduction of $J^{n+}$ and Theorem E.
\section{Arnold's $J^+$-invariant and its properties}
\begin{figure}[h]
\renewcommand*\figurename{Figure}
\begin{center}
\includegraphics[width=10cm]{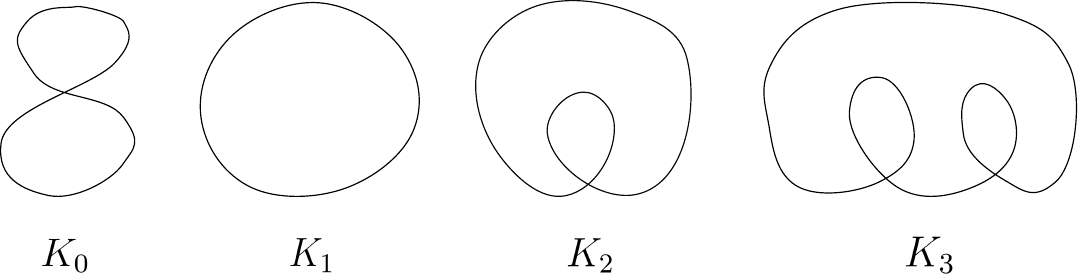}
\caption{Standard immersions $K_j$ for i = 0, 1, 2 ,3}
\label{K_j}
\end{center}
\end{figure}

\begin{figure}[h]
\renewcommand*\figurename{Figure}
\begin{center}
\includegraphics[width=12cm]{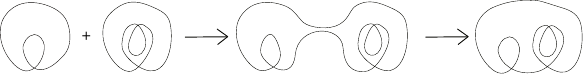}
\caption{Connected sum of $K_2$ and $L_2$}
\label{directsumex}
\end{center}
\end{figure}

\begin{figure}[h]
\renewcommand*\figurename{Figure}
\begin{center}
\includegraphics[width=9cm]{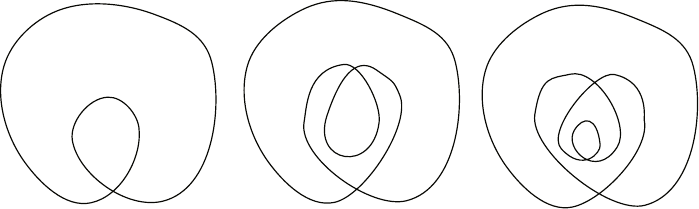}
\caption{Immersions $L_1$, $L_2$ and $L_3$}
\label{L_j}
\end{center}
\end{figure}

Consider immersions of the circle $S^1$ (or $\mathbb{R} / \mathbb{Z}$) in the plane and define these as smooth functions ${S: S^1 \rightarrow \mathbb{C}}$, whose derivative never vanishes. We use $S$ also as notation of the image of this function. In the following we examine generic immersions and an immersion is said to be generic if it has only transverse double points, which we denote as $p$.
\begin{figure}[h]
\renewcommand*\figurename{Figure}
\begin{center}
\includegraphics[width=8cm]{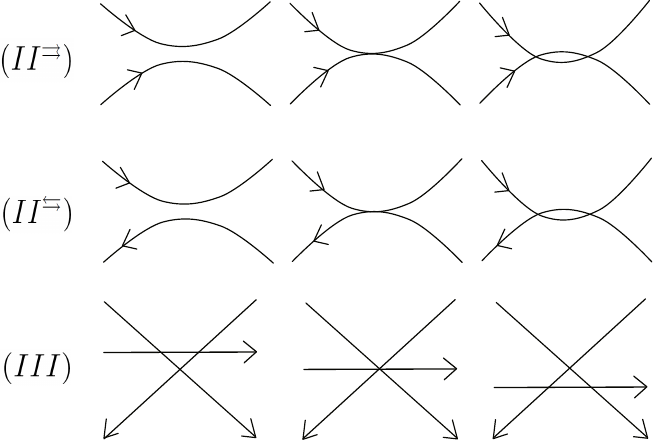}
\caption{Critical scenarios during homotopy}
\label{scenarios}
\end{center}
\end{figure}
\begin{defi}
Two immersions $S$ and $S'$ are \textbf{homotopic}, if there exists a smooth function \\
${g: {S}^{1} \times [0,1] \rightarrow \mathbb{C}}$ with $g(\cdot,0)=S$ and $g(\cdot,1)=S'$  such that for all $t \in [0,1]$ $g(\cdot,t): {S}^{1} \rightarrow \mathbb{C}$ is an immersion.
\end{defi}
During homotopies of generic immersions there are three critical homotopy scenarios: Double points can vanish or new ones can occur by passing direct ($II^{\rightrightarrows}$) or inverse ($II^{\leftrightarrows}$) self-tangencies or the immersion can pass through triple points ($III$). In [1] Arnold introduces three invariants, which are each sensible to one of these three scenarios visualized in Figure \ref{scenarios}. Positively passing direct or inverse self-tangency leads to new double points, whereas negatively passing means that double points vanish. Before and after going through a triple point, the immersion has the same amount of double points.\\
Arnold's  $J^+$-invariant is invariant under inverse self-tangency and triple points and sensible to direct self-tangency. \\
Let $S$ be a generic immersion, then
\begin{enumerate}
\item ${J}^{+}(S) \in 2\mathbb{Z}$ is well-defined and independent of orientation
\item ${J}^{+}(S)$ is constant under homotopy including triple points ($III$) and inverse self-tangency ($II^{\leftrightarrows}$)
\item ${J}^{+}(S)$ increases or decreases by the value of $2$ when traversing direct self-tangency ($II^{\rightrightarrows}$) positively or negatively
\item For predefined standard immersions $K_j$ with $j \in \mathbb{N}_0$, as visualized in Figure \ref{K_j}, $J^+$ is given through
\end{enumerate}
\[{J}^{+}(K_j)=\begin{cases}
2 - 2j & \text{for } j \neq 0\\
 0 & \text{for } j=0  .
\end{cases}\]
\ \\
Due to the Theorem of Whitney-Graustein, discussed in [7], every arbitrary immersion $S$ with winding number $j$ is homotopic to the suitable standard immersion $K_j$. Hence the $J^+$-invariant of every immersion can be calculated by finding such a homotopy transforming $S$ in $K_j$ and counting the amount of direct tangencies. 
\begin{prop} For every $x \in 2 \mathbb{Z}$ one can find an immersion $S$ such that $J^+(S)=x$.\end{prop}
\textit{Proof:} For $x\leq 0 $ one can choose the suitable standard immersion $K_j$. For $x > 0$, one can chose a suitable immersion $A_k$ as defined in Figure \ref{A_k} in chapter 3. $\hfill \Box$ \\
\ \\
In [6], Viro's formula is introduced as a new method to calculate the $J^+$-invariant of an arbitrary immersion. Using the formula it is not necessary to construct a homotopy to a standard immersion. Instead, one can directly calculate the $J^+$-invariant by determining topological properties. Let $S$ be an arbitrary immersion in the complex plane with double points $p$. The connected components in $\mathbb{C} \setminus \{S\}$ are called $C$ with winding number $\omega_C(S) \in \mathbb{Z}$. The index $ind_p(S) \in \mathbb{Z}$ of a double point $p$ is the arithmetic mean of the winding numbers of the four adjacent components $C_i$ for $1\leq i \leq 4$ at $p$, as formula given through $ind_p(S) := \frac{1}{4} \sum_{i=4}^4 \omega_{C_i}(S)$. 
\begin{prop}\textbf{(Viro's formula)}
The ${J}^{+}$-invariant of an immersion $S$ with $n$ double points is given through 
\begin{equation}
{J}^{+}(S) = 1 + n - \sum_{C \subset \mathbb{C}\setminus S} \omega_C(S)^2 + \sum_p ind_p(S)^2. \nonumber
\end{equation} 
\end{prop}
\begin{figure} [h]
\renewcommand*\figurename{Figure}
\begin{center}
\includegraphics[width=4cm]{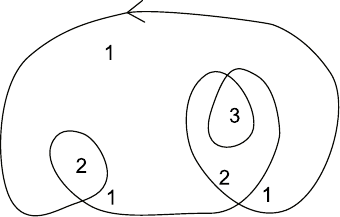}
\caption{Calculating the $J^+$-invariant using Viro's formula}
\label{Viro_ex}
\end{center}
\end{figure}
\textbf{Example:} Figure \ref{Viro_ex} shows an immersion dividing the complex plane $\mathbb{C}$ in five connected components with winding numbers $0$, $1$, $2$, $2$ and $3$. It has $n=3$ double points, two double points with the index $1$ and one with the index $2$. Using Viro's formula one calculates\\
\ \\
\begin{tabular}{ll}
${J}^{+}(S)$ & $= 1 + 3 - (0^2 + 1^2 + 2^2 + 2^2 + 3^2) + (1^2 + 1^2 + 2^2)$ \\ & $= 1 + 3 - 18 + 6$\\ & $= -8$.\end{tabular}\\
$\textit{Proof of Viro's formula: }$ We start by verifying Viro's formula for all standard immersions $K_j$. It is sufficient to take positively oriented immersions $K_j$, as Viro's formula  contains only squares of $ind_p(K)$ and $\omega_C(K)$. Then, we show that the formula behaves like Arnold's $J^+$-invariant under the critical scenarios.\\
The immersion $K_0$ divides $\mathbb{C}$ into three connected components with the winding numbers $1$, $-1$ and $0$ and has one double point $p$ $(n = 1)$ with $ind_p(K_0) = 0$. Hence, we get ${J}^{+}(K_0) = 1 + 1 - (+1)^2 - (-1)^2 + - (0)^2 + 0 = 0$.\\
Set $j=1$, $\mathbb{C} \setminus \{K_1\}$ has two components with the winding numbers $0$ and $1$ and $K_1$ has no double points $(n = 0)$. Hence, we calculate ${J}^{+}(K_1) = 1 + 0 - 1^2 - 0^2 = 0$. \\
All $K_j$ with $j > 1$ have $n = j - 1$ double points $p$ with index $ind_p(S) = 1$ and $\mathbb{C} \setminus \{K\}$ has $j + 1$ components. Hence, we get \begin{center}${{J}^{+}(K_j) = 1 + (j - 1) - (1^2 + (j - 1)\cdot 2^2 + 0^2)+ (j - 1)\cdot 1^2 = -2\cdot(j - 1)}$. \end{center}
Hence, Viro's formula for $K_j$ coincides with Arnold's definition of $J^+(K_j)$. Now let $S$ be an arbitrary immersion, which is homotopic to a standard immersion $K_j$. During a homotopy, three different critical situations can occur, which lead to a change of the topological variables that are used in Viro's formula: Crossing direct and inverse self-tangency and traversing triple points. If the immersion $S$ has just passed one of these three situations during a homotopy, we can write
\begin{align}
 J^+(S_{\text{new}}) & =  J^+(S_{\text{old}}) + \triangle {J}^{+}(S) \nonumber \\
 & = J^+(S_{\text{old}}) + \triangle n - \triangle  \sum \omega_c(S)^2 + \triangle \sum ind_p(S)^2. \nonumber
\end{align} 
\begin{figure} [h]
\renewcommand*\figurename{Figure}
\begin{center}
\includegraphics[width=4cm]{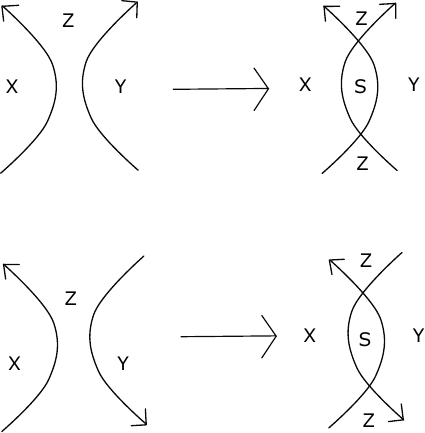}
\caption{Components when traversing self-tangency}
\label{Viro_Bew4}
\end{center}
\end{figure}\\
Now we study $\triangle J^+(S)$ for all three critical situations. Figure \ref{Viro_Bew4} shows the local neighbourhood before and after traversing direct or inverse self-tangency. In Figure \ref{Viro_Bew4}, the right and the left path can be connected in an arbitrary way. By passing a self-tangency, the immersion $S$ gets two new double points $A$ and $B$, so $\triangle n = 2$. Furthermore, two new components occur in $\mathbb{C}\setminus \{S\}$. The component $C_{s}$ is the newly arose component with winding number $s$. Besides that the component with winding number $z$, which is first in between the branches, is divided in two components with winding number $z$. Figure \ref{Viro_Bew4} illustrates the situation with winding numbers $x$, $y$, $z$ and $s$. Therefore, we have ${\triangle \sum \omega_C(S)^2 = \omega_{C_{s}}(S)^2 + \omega_{C_{z}}(S)^2}$, as these two new components occur and the winding numbers of the other components do not change. Furthermore, ${\triangle \sum ind_p(S)^2 = ind_A(S)^2 + ind_B(S)^2}$ with $ind_A(S)=ind_B(S) =\frac{1}{4}(x + y + z + s)$ follows, as $ind_p(S)$ is constant for all other double points $p$ of $S$.
\begin{figure} [h]
\renewcommand*\figurename{Figure}
\begin{center}
\includegraphics[width=4cm]{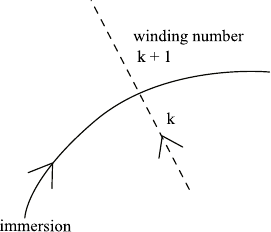}
\caption{Method to determine winding number of an arbitrary component}
\label{winding}
\end{center}
\end{figure}\\
Using the rule for determining the winding numbers of the single components, as visualized in Figure \ref{winding}, they are determinded through $x = z + 1$, $y = z - 1$ and $s = z$ for an arbitrary $z \in \mathbb{Z}$ in the case of direct self-tangency. Therefore, we conclude
\begin{align*}
\triangle {J}^{+}(S) & = 2 - (z^2  + s^2) + 2 \left(\frac{1}{4}(x + y + z + s)\right)^2  \\
 & = 2 - 2(z)^2 + 2(z)^2 = 2
\end{align*}\\
Traversing an inverse self-tangency  gives us $x = z + 1 = y$ and $s = z + 2$ for an arbitrary $z \in \mathbb{Z}$. Therefore, we get
\begin{align*}
\triangle {J}^{+}(S) & = 2 - (z^2  + s^2) + 2 \left(\frac{1}{4}(x + y
                       + z + s)\right)^2  \\
 & = 2 - (z)^2 - (z+2)^2 + 2(z+1)^2  = 0 
\end{align*}\\
\begin{figure}[h]
\renewcommand*\figurename{Figure}
\begin{center}
\includegraphics[width=14cm]{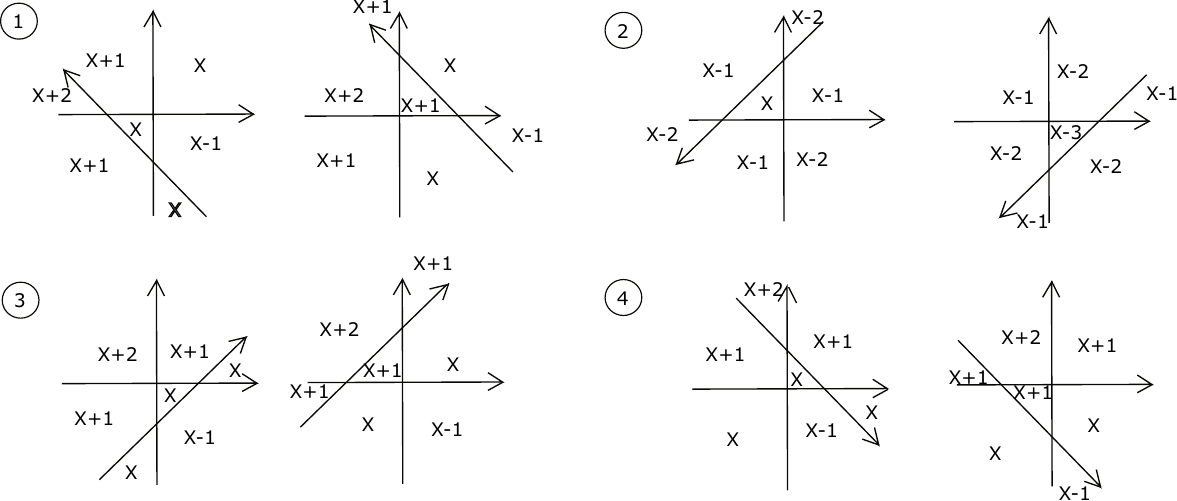}
\caption{Winding numbers when traversing triple points}
\label{Viro_Bew6}
\end{center}
\end{figure}\\
Traversing a triple point does not change the amount of double points $n$ and also the amount of components $C$ in $\mathbb{C}\setminus \{S\}$ do not change ($\triangle n =0$). Only the winding number of the new emerging triangle $T_{\mathrm{new}}$ compared to the old triangle $T_{\mathrm{old}}$ changes. Therefore, the indices of the double points $A$, $B$ and $C$, which are the corners of the triangle, also change. In Figure \ref{Viro_Bew6}, two axes are fixed orthogonally, the third axis moves through the triple point, the figure shows all four possibilities of orientations. We fix the winding number of the $T_{\mathrm{old}}$ to be $x$ and get winding numbers as given in Figure \ref{Viro_Bew6}. \\
Traversing a triple point, Viro's formula changes in the following way
\begin{align*}
\triangle {J}^{+}(S) &  = - \triangle \sum \omega_C (S)^2 + \triangle \sum ind_p(S)^2 \nonumber \\
 & = - \left( \omega_{T_{\mathrm{new}}}(S)^2 - \omega_{T_{\mathrm{old}}}(S)^2 \right) + \triangle ind_A(S)^2 + \triangle ind_B(S)^2 + \triangle ind_C(S)^2.
\end{align*}\\
\begin{tabular}{ll}
1. version: $\triangle {J}^{+}(S)$  & $ = - ((x + 1)^2 - x^2) + (x^2 + 2(x + 1)^2 - 2x^2 - (x + 1)^2 )$\\
 & $ =  x^2 - (x + 1)^2 - x^2 + (x + 1)^2 =0$\\
 & \\
2. version: $\triangle {J}^{+}(S)$&  $= - ((x - 3)^2 - x^2) + (3(x - 2)^2 - 3(x - 1)^2)$\\
 & $=  - x^2 + 6x - 9 + x^2 + (3x^2 - 12x + 12 - 3x^2 + 6x - 3)= 0$\\
  & \\
3. version: $\triangle {J}^{+}(S)$& $= - ((x + 1)^2 - x^2) + ( x^2 + 2(x + 1)^2 - 2x^2 - (x + 1)^2) $\\
 & $= - (x + 1)^2 + x^2 - x^2 + (x + 1)^2 = 0$\\
  & \\
4. version: $\triangle {J}^{+}(S)$& $= - ((x + 1)^2 - x^2) + (x^2 + 2(x + 1)^2 - 2x^2 - (x + 1)^2) $\\
 & $=  - (x + 1)^2 + x^2 - x^2 + (x + 1)^2 = 0$\\
  & \\
\end{tabular}\\
In all four versions of traversing a triple point, we get $\triangle J^+(S) =0$. Thus, Viro's formula behaves like the $J^+$-invariant. It changes by the value of $2$ through traversing direct self-tangency and is constant under traversing inverse self-tangency and triple points. $\hfill \Box$\\
\ \\
The following corollaries are easily shown using Viro's formula.
\begin{cor} The $J^+$-invariant is additive under taking connected sums.\end{cor} 
Figure \ref{directsumex} visualizes the connected sum of two immersions.
\ \\
\textbf{Example:} Figure \ref{Viro_ex} shows an immersion which is a connected sum of $K_2$(= $L_1$) and $L_2$, as visualized in Figure \ref{directsumex}. The immersions $L_j$ are defined in Figure \ref{L_j}. 
\ \\
By the definition of $L_j$ and Viro's formula, one concludes:
\begin{align*}
{J}^{+}(L_j) & = 1 + j -  \sum\nolimits_{k=1}^{j+1} k^2 + \sum\nolimits_{k=1}^{j} k^2  \nonumber \\
  & = 1 + j - (j + 1)^2  =  - j - j^2 = -j \cdot (j+i).\end{align*}  
We easily calculate that the $J^+$-invariant of the immersion in Figure \ref{directsumex} is given through \begin{center}$J^+(S) = J^+(K_2) + J^+(L_2) = (2 -2 \cdot 2) + ( -2 \cdot (1 + 2)) = -2 + ( -6) = -8$.\end{center}
\ \\
\begin{cor} The $J^+$-invariant is constant under the addition of an exterior loop.
\end{cor}
This is easily seen as the addition of an exterior loop is equivalent to adding $K_0$ with $J^+(K_0)=0$ via connected sum. This corollary will be used later when constructing pairs of immersions with suitable invariants. \\
\section{The ${J}^{2+}$-invariant for oriented pairs of generic immersions}
This chapter introduces the $J^{2+}$-invariant for oriented pairs of generic immersions $K$. This invariant is constant when traversing inverse (self-)tangency and triple points, but it changes through direct (self-) tangency (Theorem A). The relation between the single immersions is considered by the introduction of an encircling index $u(K)$ which is used to calculate the $J^{2+}$-invariant. After the introduction of the invariant Theorem A stating that the $J^{2+}$-invariant behaves as described above is proven. In section 2.2, several useful properties of the invariant will be presented, namely the independence of $J^{2+}$ and the $J^+$-invariants of the immersions forming the pair (Theorem B).
\subsection{Introduction of the $J^{2+}$-invariant}
Let $K$ be a pair of two immersions $S_1: S^1 \rightarrow \mathbb{C}$ and $S_2: S^1 \rightarrow \mathbb{C}$ defined as function
\begin{equation} K : S^1 \sqcup S^1 \rightarrow \mathbb{C}. \nonumber \end{equation}
In the following, $K$ is used as notation for the image of this function. We write $K(S_1,S_2)$ for a pair of immersions consisting of the two single immersions $S_1$ and $S_2$, whereas we have to consider that $S_1$ and $S_2$ can be paired in many different ways. Furthermore, we want $K(S_1,S_2)$ to be generic which means that it has only transverse double points. Based on Viro´s formula for $J^+$, let $n$ be the amount of double points of $K$. The double points can be double points of a single immersion $S_1$ or $S_2$ or intersection points of $S_1$ and $S_2$. The index $ind_p(K)$ of a double point $p$ and the winding number $\omega_C(K)$ of a component $C$ are defined as in Viro´s formula. To describe the relation of the immersions forming a pair towards each other, $u(K)$ is introduced.
\begin{defi} The \textbf{encircling index} of a pair $K(S_1,S_2)$ is defined by
\begin{center}
$u(K) = \begin{cases}
    0 & S_1 \cap S_2 \neq \emptyset \\
    \omega_{C_1}(S_2) + \omega_{C_2}(S_1) & \text{where } S_1 \subset C_1 \in \mathbb{C}\setminus \{S_2\} \text{ and } S_2 \subset C_2 \in \mathbb{C}\setminus \{S_1\}
\end{cases} $\end{center} \end{defi}
Note $u(K)$ is an integer and describes how often one immersion encircles the other one completely. 
\begin{defi} A pair of immersions $K(S_1,S_2)$ is called \textbf{disjoint}, if $S_1$ is contained in the unbounded component of $\mathbb{C} \setminus \{S_2\}$. \end{defi}
Two pairs of immersions $K$ and $K'$ are called \textbf{homotopic}, if there exists a smooth function $g: (S^1 \sqcup S^1) \times [0,1] \rightarrow \mathbb{C}$ with $g(\cdot,0) = K$ and $g(\cdot,1)=K'$ such that for all $t \in [0,1]$ $g(\cdot,t): S^1 \sqcup S^1 \rightarrow \mathbb{C}$ is a generic pair of immersions except at finitely many $s \in (0, 1)$ at which the three critical scenarios can occur. \\
We call a homotopy a \textbf{single-homotopy} of a pair of immersions $K(S_1,S_2)$, if the critical scenarios only include double points of $S_1$ (or $S_2$). Hence, the relation between the two immersions does not change ($\Delta u (K)=0$). We call a homotopy a \textbf{unravel-homotopy} of a pair of immersions $K(S_1,S_2)$, if the critical scenarios only include double points between $S_1$ and $S_2$. It is obvious, that every arbitrary pair of immersions $K(S_1, S_2)$ is homotopic to a disjoint pair of immersions $K'(S_1,S_2)$ via unravel-homotopy and that every disjoint pair $K'(S_1,S_2)$ is homotopic to a disjoint pair out of standard immersions $K'(K_i, K_j)$ for some $i,j \in \mathbb{N}_0$ via single-homotopy.
\begin{defi} The $J^{2+}$-invariant for an oriented pair of immersions $K$ is defined by
\begin{equation}
J^{2+}(K) := 2 + n - \sum_{C\subset \mathbb{C}\setminus K} \omega_c(K)^2 + \sum_p ind_p(K)^2 + u(K)^2. \nonumber
\end{equation} \end{defi}
\begin{cor} The $J^{2+}$-invariant of a disjoint pair of immersions $K(S_1,S_2)$ is the sum of the $J^+$-invariants of $S_1$ and $S_2$:
\begin{equation}{J}^{2+}(K) = {J}^{+}(S_1) + {J}^{+}(S_2). \nonumber \end{equation}
\end{cor}
\textit{Proof: } The encircling index of a disjoint pair is $u(K)=0$ and by definition of $J^{2+}$, the formula can be split into Viro's formulas for $J^+(S_1)$ and $J^+(S_2)$. $\hfill \Box$\\
\begin{theo}[Theorem A] The $J^{2+}$-invariant is constant under inverse tangency and triple points and increases (decreases) by the value of $2$ when going positively (negatively) through direct tangencies. \end{theo}
\begin{figure} [h]
\renewcommand*\figurename{Figure}
\begin{center}
\includegraphics[width=13cm]{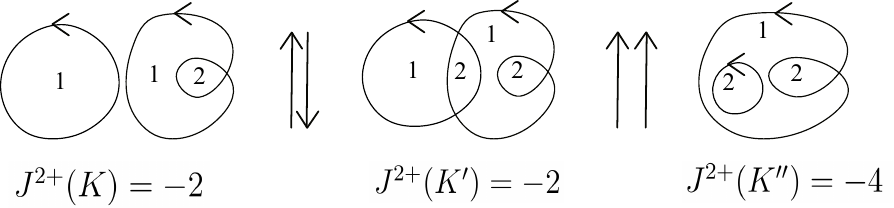}
\caption{Change of $J^{2+}$-invariant during homotopy of a pair of immersions}
\label{J2+Ex}
\end{center}
\end{figure}
\textbf{Example:} Figure \ref{J2+Ex} visualizes an example of a homotopy of a pair of immersions. The pair of immersions $K$ is first transformed into $K'$ by pushing the single immersions into each other going through positive inverse and negative direct tangency until one immersions encircles the other one completely in $K''$. The $J^{2+}$-invariant does not change by positive inverse tangency and changes by the value of $-2$ when going through negative direct tangency.
\begin{center} \begin{tabular}{ll}$J^{2+}(K)$ & $=2 + 1 - ( 1^2 + 1^2 + 2^2) + 1^2 + 0^2 = -2$\\
$J^{2+}(K')$ & $=2 + 3 - ( 1^2 + 1^2 + 2^2+2^2) + (1^2+ 1^2+ 1^2) + 0^2 = -2$\\
$J^{2+}(K'')$ & $=2 + 1 - ( 1^2 + 2^2+2^2) + 1^2 + 1^2 = -4$\end{tabular}\end{center}

\ \\
\textit{Proof of  Theorem A:} It is obvious, that the $J^{2+}$-invariant is well-defined for pairs of immersions $K(S_1,S_2)$. In the following we show that $J^{2+}$ behaves correctly under critical scenarios. Lemma 2.6 discusses the behavior of $J^{2+}$ during single-homotopy. The next two Lemmata discuss the behavior under critical scenarios during unravel-homotopy. Lemma 2.7 discusses tangencies between $S_1$ and $S_2$ and Lemma 2.8 studies the behavior when going through triple points. 
\begin{lem} Under single-homotopy, the $J^{2+}$-invariant of an oriented pair $K(S_1,S_2)$ changes by the value of $2$ when going through direct tangency ($II^\rightrightarrows$) and remains constant when traversing inverse tangency ($II^\leftrightarrows$) or triple points ($III$).\end{lem}
\textit{Proof:} The Proof of Viro's formula can be applied analogously. We know $\Delta u(K) =0$, as a single-homotopy does not change the relation between $S_1$ and $S_2$. Thus, we can locally apply the proof of Viro's formula.$\hfill \Box$\\
\ \\
The next two lemmata study the bahavior of $J^{2+}$ under the critical scenarios within a unravel-homotopy. In such cases the encircling index $u(K)$ might change.
\begin{lem} During an unravel-homotopy, the $J^{2+}$-invariant changes by the value of $2$ when going through direct tangency ($II^\rightrightarrows$) and stays constant when traversing inverse tangency ($II^\leftrightarrows$).\end{lem}
\begin{figure} [h]
\renewcommand*\figurename{Figure}
\begin{center}
\includegraphics[width=4cm]{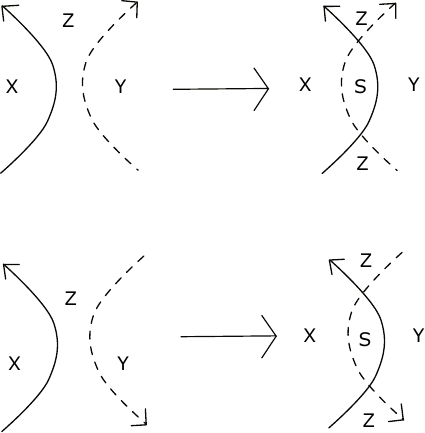}
\caption{Local neighbourhood when traversing a tangency of $S_1$ and $S_2$}
\label{allgFall}
\end{center}
\end{figure}
\textit{Proof: }Let $K$ be an oriented pair of immersions, which traverses a tangency between $S_1$ and $S_2$ during a homotopy. Figure \ref{allgFall} shows the local neighbourhood by traversing a direct or inverse tangency. The left branch is part of the immersion $S_1$, the right one is part of the immersion $S_2$. The values $x$, $y$, $z$ and $s$ are the winding numbers of the single components. We consider
\begin{equation*}
\Delta J^{2+}(K) = \Delta n - \Delta \sum_{C \subset \mathbb{C} \setminus K} {\omega_C(K)}^2 + \Delta \sum_p {ind_p(K)}^2 + \Delta u(K)^2.
\end{equation*}
By positively traversing a tangency we get two new double points $p$ with $ind_p(K) = \frac{1}{4} (x + y + z + s)$, hence we obtain
$\Delta n = 2$ and $\Delta \sum_p ind_p(K) ^2 = 2 (\frac{1}{4}(x + y + z + s))^2$, as the indices of all other double points of $K$ do not change. Considering $\Delta \sum_{C \subset \mathbb{C} \setminus K} \omega_C(K) ^2$ and $\Delta u(K)$ we need to distinguish two cases:
\begin{itemize}
    \item[1.] If $S_1\cap S_2 \neq \emptyset$ before the move, $u(K)$ does not change ($\Delta u(K)^2=0$) and there are two more new connected components in $\mathbb{C}\setminus K$ ($\Delta \sum_{C \subset \mathbb{C} \setminus K} \omega_C(K)^2 =s^2+z^2$) 
    \item[2.] If $S_1\cap S_2 = \emptyset$ before the move, $u(K)$ becomes $0$ ($\Delta u(K)^2 = - z^2$)  and there is one new connected component in $\mathbb{C}\setminus K$ ($\Delta \sum_{C \subset \mathbb{C} \setminus K} \omega_C(K) ^2 = s^2$)
\end{itemize}
Therefore, $\Delta J^{2+}$ is given through
\begin{center}
    $\Delta J^{2+}(K) = 2 + \frac{1}{8}(x + y+ z +s)^2-s^2-z^2$
\end{center}
For a fixed $z \in \mathbb{Z}$ we determine for the winding numbers $x = z + 1$, $y = z - 1$ and $s=z$ for traversing a direct tangency, whereas for traversing an inverse tangency we know that $x = z + 1 = y$ and $s = z + 2$. Therefore, we receive
\begin{center}$\Delta {J}^{2+}(K) = 2 +  \frac{1}{8}(4z)^2 - z^2 - z^2= 2$\end{center} 
for positively traversing a direct tangency and 
\begin{center}$\Delta {J}^{2+}(K) = 2 + \frac{1}{8}(4z + 4)^2 - (z+2)^2 - z^2 = 0$\end{center} 
for positively traversing an inverse tangency. For negatively traversing direct or inverse tangency $\Delta J^{2+}(K)=-2$ and $\Delta J^{2+}(K)=0$ follow analogously.$\hfill\Box$\\
\begin{lem} The ${J}^{2+}$-invariant is constant when going through triple points of $S_1$ and $S_2$.\end{lem}
\textit{Proof:} The immersions $S_1$ and $S_2$ intersect each other before and after the move, hence, $\Delta u(K)^2=0$. Thus the proof works similar as for $J^+$, see section 1.2 and Figure \ref{Viro_Bew6}. $\hfill \Box$\\
\ \\
Combining the three lemmata, it is shown that the $J^{2+}$-invariant increases (decreases) by the value of $2$ when going through direct tangency and remains constant when traversing inverse tangency and triple points. $\hfill\Box$\\
\subsection{Properties of the ${J}^{2+}$-invariant}
\begin{figure}[h]
\renewcommand*\figurename{Figure}
\begin{center}
\includegraphics[width=16cm]{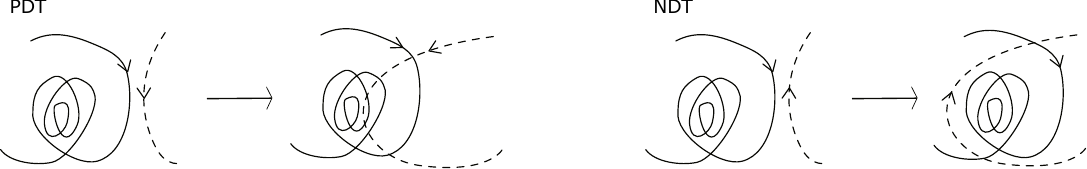}
\caption{Moves to create $j$ positive or negative direct tangencies}
\label{PDT}
\end{center}
\end{figure}

This chapter highlights several properties of the $J^{2+}$-invariant. Especially interesting is Theorem B which states that the $J^{2+}$-invariant of an oriented pair $K(S_1,S_2)$ is independent of the $J^+$-invariants of the single immersions. The proof of the theorem is an algorithm constructing a suitable pair of immersions for given $J^{2+}$- and $J^+$-invariants.
\begin{cor} The $J^{2+}$-invariant is independent of orientation for disjoint pairs of immersions. \end{cor}
\textit{Proof:} This is true due to the above proposition and the orientation independence of Arnold's $J^{2+}$-invariant.$\hfill \Box$\\
\begin{prop}
For all arbitrary pairs of immersions $K$, we have ${J}^{2+}(K) \in 2\mathbb{Z}$.
\end{prop}
\textit{Proof:} This is trivial, as the $J^{2+}$-invariant only changes by the value of $2$ and as every oriented pair $K$ is homotopic to a disjoint pair $K'$ with a $J^{2+}$-invariant given through the sum of the $J^+$-invariants of the single immersions which are in $2\mathbb{Z}$. $\hfill\Box$\\
In comparison to the $J^+$-invariant, the $J^{2+}$-invariant is not independent of orientation, but the following proposition holds. 
\begin{prop} The $J^{2+}$-invariant is independent under simultaneous change of orientation of both immersions. \end{prop}
An example for an simultaneous orientation change is visualized in Figure \ref{simultan} in the next chapter. \\
\textit{Proof:} Let $K$ be an arbitrary oriented pair of immersions, and we change the orientation of both immersions simultaneously. The amount of double points and $u(K)$ do not change. The winding numbers of all components, as well as the indices of the double points change their sign. But as in the formula to calculate $J^{2+}$ all these values are summed squared, $J^{2+}(K)$ does not change.$\hfill \Box$\\
\begin{prop} The $J^{2+}$-invariant is constant under the addition of an exterior loop.\end{prop}
\textit{Proof:} Adding an exterior loop is adding $K_0$ via connected sum at an unbounded component. The pair gets one new double point with index $0$ and an additional component with winding number $\pm 1$, $u(K)$ does not change and, thus, one receives
\begin{center}\begin{tabular}{ll} $\Delta J^{2} (K)$ & $= \Delta n - \Delta \sum_{C \subset \mathbb{C} \setminus K} \omega_C(K)^2 + \Delta \sum_p ind_p(K)^2 + \Delta u(K)^2$ \\
 & $= 1 - (\pm1)^2 + 0^2 + 0^2$ \\ & $= 0$. \end{tabular} \end{center} $\hfill \Box$\\
The $J^+$-invariants of the single immersions forming a pair are independent of the $J^{2+}$-invariant.
\begin{theo}[Theorem B] For arbitrary even numbers $x$, $y$ and $z$ there always exists an oriented pair of generic immersions $K(S_1,S_2)$ such that ${J}^{+}(S_1)=x$, ${J}^{+}(S_2)=y$ and ${J}^{2+}(K)=z$. \end{theo}
\textit{Proof:} The proof is given through \textit{Algorithm 0} which constructs suitable pairs $K(S_1,S_2)$ for arbitrary variables $x$, $y$ and $z$ such that $J^+(S_1) = x$, $J^+(S_2)=y$ and $J^{2+}(K)=z$. \\
\ \\
We take two immersions $S_1$ and $S_2$ with suitable $J^+$-invariant, whereas $S_2$ can be chosen arbitrarily and $S_1$ contains some part $L_j$, whereas $L_j$ is defined as in Figure \ref{L_j}. $S_1$ and $S_2$ form a disjoint pair $K'(S_1,S_2)$ with $J^{2+}(K') = x + y$. When pushing the outer branch of a second immersion into the $L_j$-part one creates positive and negative direct tangencies depending on the orientation of the single immersions. These moves are called $PDT$ and $NDT$ and they are visualized in Figure \ref{PDT}. Denote $K$ to be the resulting pair after the moves.\\
\begin{itemize}
\item [\textbf{PDT}] positive direct tangencies:\\
The two immersions $S_1$ and $S_2$ face each other in components being differently oriented. The immersion $S_1$ contains the immersion $L_j$. The outer branch of $S_2$ is pushed into $L_j$ to create $j$ positive direct tangencies, thus the $J^{2+}$-invariant changes by the value $2j$. 
\item [\textbf{NDT}] negative direct tangencies:\\
The two immersions $S_1$ and $S_2$ face each other in components being equally oriented. The immersion $S_1$ contains the immersion $L_j$.  The outer branch of $S_2$ is pushed into $L_j$ to create $j$ negative direct tangencies, to do so the pair has to go through $j+1$ positive inverse tangencies. The $J^{2+}$-invariant changes by the value $-2j$.
\end{itemize}
The following algorithm is a version of the algorithm given in [4] to construct suitable pairs for arbitrary $x$, $y$ and $z \in 2 \mathbb{Z}$ such that $J^+(S_1) = x$, $J^+(S_2)=y$ and $J^{2+}(K)=z$.\\
\ \\
\hrule
\ \\
\textbf{Algorithm 0:}\\
For given arbitrary even numbers $x$, $y$ and $z \in 2\mathbb{Z}$ calculate $d:= z - x - y$.
\begin{itemize}
\item[1.] If $d = 0$, which is equivalent to $z = x + y$, construct a disjoint pair of immersions out of $S_1$ with $J^+(S_1) = x$ and $S_2$ with $J^+(S_2) = y$. One calculates \\
\begin{center}$J^{2+}(K) = J^+(S_1) + J^+(S_2) = x + y = z$.\end{center} 
\item[2.] If $d \neq 0$: \\
Define the immersion $S_1$ as a connected sum of the immersions $S_{1,1}$ and $S_{1,2}$, where $S_{1,2}$ is $L_{\frac{|d|}{2}}$ and $S_{1,1}$ is defined as immersions such that 
\begin{equation}{J}^{+}(S_{1,1}) = x - {J}^{+}\left(L_{\frac{|d|}{2}}\right) = x - (- \left(\frac{|d|}{2}\right) - \left(\frac{|d|}{2}\right)^2).\nonumber\end{equation}
Such an immersion $S_{1,1}$ exists as $x, {J}^{+}\left(L_{\frac{|d|}{2}}\right) \in 2\mathbb{Z}$ and due to Proposition 1.2. Hence, 
\begin{equation}J^+(S_1) = {J}^{+}(S_{1,1}) + {J}^{+}(S_{1,2}) = x. \nonumber\end{equation}
Next choose an arbitrary immersion $S_2$ with $J^+(S_2) = y$.\\
Let $K'$ be a pair of disjoint immersions $S_1$ and $S_2$, thus the $J^{2+}$-invariant is given through \begin{equation}J^{2+} (K')=x + y.\nonumber\end{equation}
Choose the orientations of $S_1$ and $S_2$ according to the moves \textbf{PDT} and \textbf{NDT}.
\item[2.1] For $d > 0$, use move \textbf{PDT}. The outer branch of $S_2$ is pushed into $S_{1,2}$ such that $K$ traverses  $\frac{|d|}{2}$ direct tangencies positively, therefore the ${J}^{2+}$-invariant increases by the value of $2 \cdot \frac{|d|}{2} = |d| = d$. For the constructed system $K$ one calculates
\begin{equation}{J}^{2+}(K) ={J}^{2+}(K') + \triangle J^{2+}(K) = x + y + d = z. \nonumber\end{equation}
\item[2.2] For $d < 0$, use move \textbf{NDT}. The outer branch of $S_2$ pushes into $S_1$ until $\frac{|d|}{2} + 1$ inverse tangencies are traversed positively, $\frac{|d|}{2}$ triple points points are traversed and $\frac{|d|}{2}$ direct tangencies are traversed negatively. Hence, the ${J}^{2+}$-invariant decreases by the value of $ 2 \cdot \frac{|d|}{2} = |d| = -d$, therefore conclude 
\begin{equation}{J}^{2+}(K) = {J}^{2+}(K') + \triangle J^{2+}(K) = x + y - (-d) = z. \nonumber\end{equation}
\end{itemize} 
\hrule 
\ \\
\textit{Algorithm 0} proves the Theorem. $\empty \hfill \Box$\\
\section{Tuple of $J^{2+}$-invariants for pairs of immersions}
As explained in the previous chapter, the $J^+$-invariant is independent of the choice of orientation, whereas the $J^{2+}$-invariant is not, but it remains constant under a simultaneous change of orientation of both immersions. Hence, one can define two $J^{2+}$-invariants for each pair of immersions depending on the orientation of the single immersions.
\begin{figure}[h]
\renewcommand*\figurename{Figure}
\begin{center}
\includegraphics[width=15cm]{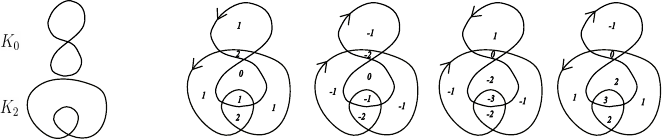}
\caption{Pairs of immersions with different orientations}
\label{simultan}
\end{center}
\end{figure}
Figure \ref{simultan} shows pairs of immersions $K$ out of $K_0$ and $K_2$, which have the $J^+$-invariants $J^+(K_0)=0$ and $J^+(K_2)=-2$ independent of their orientation. Pushing the immersions into each other, the parts pushed into each other in the first two pairs have opposite orientation, thus, there are new components with absolute winding numbers 0, 1 and 2. In the third and fourth pair, the intersecting components are equally oriented thus there are components with absolute winding numbers 0, 2 and 3. When calculating the $J^{2+}$-invariant the winding numbers of components and the indices of double points are squared. Thus, the first and second as well as the third and fourth pair have the same $J^{2+}$-invariant. Let us call the first and the second pair $K$, the third and fourth pair are denoted as $\overline{K}$.  The tuple of $J^{2+}$-invariants is given through $J^{2+}(K) = 2$ and $J^{2+}(\overline{K})=-2$. 
\begin{defi} Let $K$ be an oriented pair of immersions, we define \textbf{$\overline{K}$} to be $K$ with changed orientation of one immersion.\end{defi}
\textbf{Remark:} As the $J^{2+}$-invariant is invariant under simultaneous change of orientation, it is insignificant which immersion in $\overline{K}$ is oriented differently in comparison to the orientations in $K$. If two immersions are oriented similarly or not is not well-defined. For example $K_0$ can not be uniquely described as positively or negatively oriented. 

\subsection{Dependency of the two $J^{2+}$-invariants} 
In the following we study the dependence of the $J^{2+}(K)$ and $J^{2+}(\overline{K})$ for an oriented pair of immersions $K(S_1,S_2)$.
\begin{theo}[Theorem C] Let $K(S_1,S_2)$ be an arbitrary oriented pair, then the invariants  $J^+(S_1)=x$, $J^+(S_2)=y$, $J^{2+}(K)=z$ and $J^{2+}(\overline{K})=\overline{z}$ are even numbers satisfying $z + \overline{z} \geq 2 x + 2 y.$ \end{theo}
\textit{Proof:} It follows a proof by contradiction. Take an arbitrary pair of immersions $K(S_1, S_2)$ with even numbers $x$, $y$, $z$ and $\overline{z}$ as invariants $J^+(S_1)=x$, $J^+(S_2)=y$, $J^{2+}(K)=z$ and $J^{2+}(\overline{K})=\overline{z}$. Assume without loss of generality 
\begin{center}$z \leq x + y$, thus $z = x + y -d$ for some $d \in 2\mathbb{N}_0:= \mathbb{N} \cup \{0\}$\end{center}  and \begin{center}$\overline{z} < 2x + 2y - z$, thus $\overline{z} =  x+y + d - d^*$ for some $d^* \in 2\mathbb{N}$.\end{center} 
Every pair of immersions $K$ out of immersions $S_1$ and $S_2$ is transformable into a disjoint pair $K'$ out of the immersions $S_1$ and $S_2$ by homotopy. Then the disjoint pair $K'$ has the invariants
\begin{center}$J^{2+}(K') = J^{2+}(\overline{K'})=x + y$.\end{center}
During the homotopy from $K$ to $K'$ the amount of positive and negative direct tangencies are referred to as $t_{pd}$ and $t_{nd}$, whereas $t_{pi}$ and $t_{ni}$ are the amount of positive and negative inverse tangencies. When considering $\overline{K}$, the pair of immersions with one immersion oriented differently, the amount of tangencies during the homotopy are analogously called $\overline{t_{pd}}$, $\overline{t_{nd}}$, $\overline{t_{pi}}$ and $\overline{t_{ni}}$. We choose the homotopy to be an unravel-homotopy, thus, there is no change in the topological properties of the single immersions $S_1$ and $S_2$, only their relation towards each other changes. All the tangencies during the homotopy are tangencies between the two immersions $S_1$ and $S_2$, hence, all arising and vanishing double points are double points between the immersions. \\
Direct tangencies during the homotopy of $K$  are inverse tangencies in the homotopy of $\overline{K}$, thus the following equations hold
\begin{center}\begin{tabular}{ll} $t_{pd} = \overline{t_{pi}}$ & \\ $t_{nd}= \overline{t_{ni}}$       &    (*)\\ $t_{pi}=\overline{t_{pd}}$ & \\ $t_{ni}=\overline{t_{nd}}$. &  \end{tabular} \end{center}
During the homotopy of $K$ to $K'$ the $J^{2+}$-invariant changes from $x+y-d$ to $x+y$, whereas $d \in 2\mathbb{N}_0$. The following formula holds: 
\begin{center}\begin{tabular}{lll}
 $J^{2+}(K') $ & $= J^{2+}(K) + 2 \cdot t_{pd} - 2 \cdot t_{nd} + 0 \cdot (t_{pi} + t_{ni})$ &  \\
 & $=x + y - d + 2 \cdot t_{pd} - 2 \cdot t_{nd}$ & $\overset{!}{=}  x  + y$.\end{tabular}
\end{center}
This leads to $t_{pd} = \frac{d}{2} + t_{nd}$ which is by (*) equal to $\overline{t_{pi}} = \frac{d}{2} + \overline{t_{ni}}$.\\
When studying the homotopy of $\overline{K}$, the $J^{2+}$-invariant changes from $x+y+d+d^*$ to $x+y$,
\begin{center}\begin{tabular}{lll}
 $J^{2+}(\overline{K'}) $ & $= J^{2+}(\overline{K}) + 2 \cdot \overline{t_{pd}} - 2 \cdot \overline{t_{nd}} + 0 \cdot (\overline{t_{pi}} + \overline{t_{ni}})$ & \\
 & $=x+y + d - d^* + 2 \cdot \overline{t_{pd}} - 2 \cdot \overline{t_{nd}}$ & $\overset{!}{=}  x  + y$ \end{tabular}
\end{center}
which leads to $\overline{t_{pd}} = \frac{d^* - d}{2} + \overline{t_{nd}}$ which is with (*) equal to $t_{pi} = \frac{d^* - d}{2} + t_{ni}$.\\
In the resulting disjoint pair $K'$ the immersions do not intersect each other, hence all the double points between the immersions have to vanish. While going through positive inverse and direct tangencies the pair gains two double points between the immersions, these double points and the already existing ones must all vanish through negative inverse and direct tangencies, therefore the following estimate holds
\begin{center}
$t_{pd}+t_{pi} \leq t_{nd} + t_{ni}$ and analogously 
$\overline{t_{pd}}+\overline{t_{pi}} \leq \overline{t_{nd}} + \overline{t_{ni}}$. \end{center}
Putting the results from above into the inequality restriction, one receives
\begin{center} $t_{pd} + t_{pi} = \frac{d}{2} + t_{nd} + \frac{d^* - d}{2} + t_{ni} = t_{nd} + t_{ni} + \frac{d^*}{2} \leq t_{nd} + t_{ni}$\end{center}
and analogously\begin{center}
$\overline{t_{pd}} + \overline{t_{pi}} = \frac{d^* - d}{2} +  \overline{t_{nd}} + \frac{d}{2} + \overline{t_{ni}} = \overline{t_{nd}} + \overline{t_{ni}} + \frac{d^*}{2} \leq \overline{t_{nd}} + \overline{t_{ni}}$,\end{center}
which leads to a contradiction for $d^* > 0$ and proves the theorem.  $\hfill \Box$\\
\ \\
The following section introduces Algorithms to construct suitable pairs for given invariants satisfying $z+\overline{z} \geq 2x+2y$.

\subsection{Algorithms to construct pairs of immersions with suitable $J^{2+}$-invariants}
\begin{figure}[h]
\renewcommand*\figurename{Figure}
\begin{center}
\includegraphics[width=7cm]{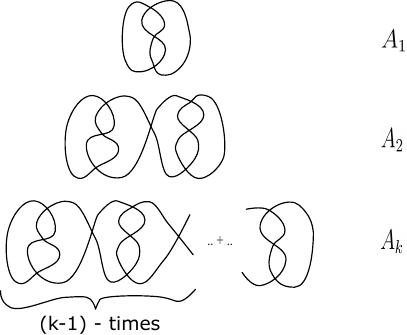}
\caption{Definition of the immersion $A_k$}
\label{A_k}
\end{center}
\end{figure}

Before introducing a new Algorithm we study the output of \textit{Algorithm 0} regarding to the invariants $J^{2+}(K)$ and $J^{2+}(\overline{K})$. Without loss of generality we define $K$ to be oriented such that the two parts facing each other and then being pushed into each other during the Algorithm to be equally oriented like visualized in \textbf{NDT}. \\
So if $d<0$ using move \textbf{NDT} the pair of immersions goes through several critical events. When calculating $J^{2+}(K)$, the pairs goes positively through inverse self tangency and triple points before going negatively through direct self tangency. Changing the orientation and doing the same moves all inverse tangencies become direct tangencies and all direct tangencies become inverse tangencies. As it is necessary to go through $\frac{|d|}{2}+1$ inverse tangencies positively for being able to go through $\frac{|d|}{2}$ direct tangencies negatively. The invariants are given through \\
\begin{center}\begin{tabular}{lll}
$J^{2+}(K') \rightarrow$ & $J^{2+}(K') - |d| $ & $ = J^{2+}(K)$\\
$J^{2+}(\overline{K'}) \rightarrow$ & $J^{2+}(\overline{K'}) +|d| +2        $ & $ = J^{2+}(\overline{K})$.
\end{tabular}\end{center}
In the case $d>0$, which uses move \textbf{PDT}, the algorithm creates positive direct tangencies in the homotopy of $\overline{K}$ without going through any other critical scenario. When considering the same moves with different orientation, namely $K$, the direct tangency is an inverse tangency, thus, $J^{2+}(K)$ remains constant \\
\begin{center}\begin{tabular}{lll}
$J^{2+}(\overline{K'}) \rightarrow$ & $J^{2+}(\overline{K'}) + d   $ & $ = J^{2+}(\overline{K})$\\
$J^{2+}(K') \rightarrow$ & $J^{2+}(K)         $ & $ = J^{2+}(K)$.
\end{tabular}\end{center}
Hence, \textit{Algorithm 0} constructs for arbitrary $x$, $y$ and $z$  $\in 2\mathbb{Z}$ and $d := z - x - y$ suitable pairs of immersions with invariants: \\
\begin{center}
\begin{tabular}{l|c|c|c|c}
 & $J^{+}(S_1)$ & $J^+(S_2)$ & $J^{2+}(K)$ & $J^{2+}(\overline{K})$\\ \hline
d=0 & $x$ & $y$ & $z$ & $z$ \\ \hline
d>0 & $x$ & $y$ & $x + y$ & $z$ \\ \hline
d<0 & $x$ & $y$ & $z$ & $2x + 2y - z + 2$ \end{tabular}\end{center}
\ \\
The table above shows that move \textbf{PDT} changes $J^{2+}(\overline{K})$ without changing the $J^{2+}(K)$, whereas move \textbf{NDT} changes both invariants. \\
\ \\
Combining the moves \textbf{PDT} and \textbf{NDT} one can only construct suitable pairs of immersions with 
\begin{center}$J^{2+}(K) \leq x + y$ and $J^{2+}(\overline{K}) \geq 2x + 2y - J^{2+}(K) + 2$\end{center}
by first constructing a pair with suitable $J^{2+}(K)$ and then adding the necessary amount of direct tangencies with move \textbf{PDT} to get the right $J^{2+}(\overline{K})$. \\
\ \\
With the introduction of new construction moves \textbf{M1, M2, M3} and \textbf{M4} one can find pairs of immersions for more combinations of invariants $z$ and $\overline{z}$. For arbitrary $x$, $y$, $z$ and $\overline{z}$ $\in \mathbb{Z}$, set $d:= z - x - y$ and $\overline{d}:= \overline{z} - x - y$.\\
\begin{figure}[h]
\renewcommand*\figurename{Figure}
\begin{center}
\includegraphics[width=16cm]{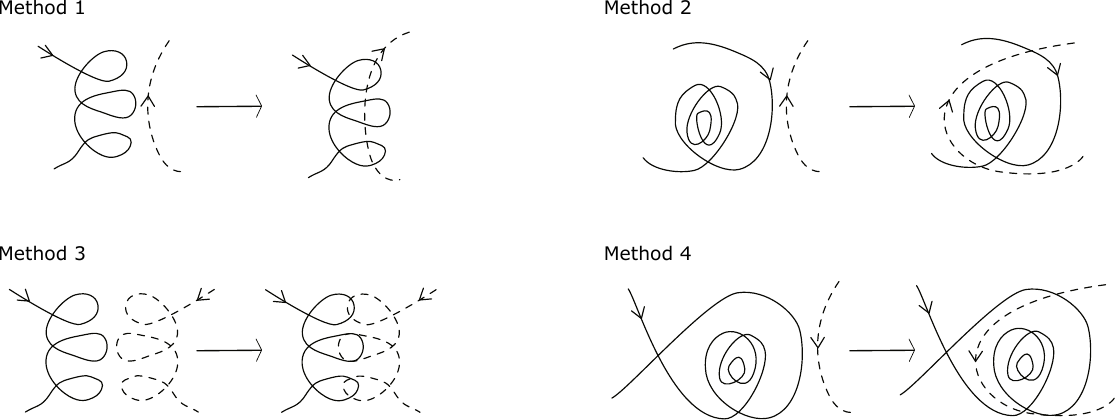}
\caption{Moves 1 to 4 with equally or differently oriented immersions $S_1$ and $S_2$ to create positive or negative direct tangencies}
\label{M1-4}
\end{center}
\end{figure}
\begin{itemize}
\item [\textbf{M1}] $d = 2k$ for some $k \in \mathbb{N}$ (equivalent to $z \geq x + y$):\\
Add $k$ exterior loops to the immersion $S_1$ and push $S_2$ into the those to create $k$ positive direct tangencies. Then $J^{2+}(K)$ changes by the value of $2k$. In comparison $\overline{K}$ goes through $k$ inverse tangencies, thus $J^{2+}(\overline{K})$ does not change.  \begin{center} \begin{tabular} {c|c}  $\Delta J^{2+}(K)$ & $\Delta J^{2+}(\overline{K})$ \\ \hline 2k & 0\end{tabular} \end{center}
\item [\textbf{M2}] $d = -2k$ for some $k \in \mathbb{N}$ (equivalent to $z \leq x + y$):\\
Use the already known construction such that one part of immersion $S_1$ is a suitable $L_k$ and push $S_2$ into $L_k$ to create $k$ negative direct tangencies. Then $J^{2+}(K)$ changes by the value of $-2k$. In comparison $\overline{K}$ goes through $k+1$ direct tangencies, thus $J^{2+}(\overline{K})$ changes by $2(k+1)$.  \begin{center}\begin{tabular} {c|c} $\Delta J^{2+}(K)$ & $\Delta J^{2+}(\overline{K})$ \\ \hline -2k & 2k + 2 \end{tabular}\end{center}
\item [\textbf{M3}] $\overline{d} = 2k$ for some $k \in \mathbb{N}$ (equivalent to $\overline{z} \geq x + y$):\\
Consider $\overline{K}$. Add $k$ exterior loops to both immersions $S_1$ and $S_2$ and push them into each other to create $k$ positive direct tangencies. Then $J^{2+}(\overline{K})$ changes by the value of $2k$. In comparison $K$ goes through $k$ inverse tangencies, thus $J^{2+}(K)$ does not change.  \begin{center} \begin{tabular} {c|c} $\Delta J^{2+}(K)$ & $\Delta J^{2+}(\overline{K})$ \\ \hline $0$ & $2k$ \end{tabular} \end{center}
\item [\textbf{M4}] $\overline{d} = -2k$ for some $k \in \mathbb{N}$ (equivalent to $\overline{z} \leq x + y$):\\
Consider $\overline{K}$. Construct $S_1$ such that it has an exterior loop including an inner $L_k$ as shown in Figure \ref{M1-4} and push $S_2$ into this $L_k$ in the exterior loop to create $k$ negative direct tangencies. Then $J^{2+}(\overline{K})$ changes by the value of $-2k$. In comparison $K$ goes through $k+1$ direct tangencies, thus $J^{2+}(K)$ changes by $2(k+1)$.   \begin{center} \begin{tabular} {c|c} $\Delta J^{2+}(K)$ & $\Delta J^{2+}(\overline{K})$ \\ \hline 2k+2 & -2k \end{tabular} \end{center}
\end{itemize}
Using the above moves one can construct pairs of immersions with an arbitrary invariants $J^{2+}(K)$ or $J^{2+}(\overline{K})$ (Theorem B). Now, we want to construct pairs $K$ with suitable $J^{2+}(K)$ and $J^{2+}(\overline{K})$. The moves above enable a construction of pairs with invariants following the restrictions given in the following theorem.\\
\begin{theo}[(Theorem D)] For arbitrary even numbers $x$, $y$ and $z$ and $\overline{z}$ satisfying $z + \overline{z} \geq 2 x + 2 y $ there always exists an oriented pair of generic immersions $K(S_1,S_2)$ such that ${J}^{+}(S_1)=x$, ${J}^{+}(S_2)=y$ and ${J}^{2+}(K)=z$ and ${J}^{2+}(\overline{K})=\overline{z}$. \end{theo}
\ \\
The theorem is proven by giving two algorithms constructing suitable pairs of immersions. \textit{Algorithm 1} uses moves \textbf{M1}, \textbf{M2}, \textbf{M3} and \textbf{M4} but excludes the cases
\begin{center}
$z < x + y$ for and $\overline{z} = 2x+2y - z$ \& $\overline{z} < x + y$ for and $z = 2x+2y - \overline{z}$.
\end{center}
These specific cases can be constructed using \textit{Algorithm 2} which uses other moves and is hence introduced separately.\\
\ \\
\hrule 
\ \\
\textbf{Algorithm 1:}\\
Let $x$, $y$, $z$ and $\overline{z}$ be arbitrary even numbers with
\begin{itemize}
    \item[(i)]$z, \overline{z} \geq x + y$, 
    \item[(ii)]$z < x + y$ and $\overline{z} > 2x+2y - z$ or
    \item[(iii)]$\overline{z} < x + y$ and $z > 2x+2y - \overline{z}$.
\end{itemize}
Case (i): $z, \overline{z} \geq x + y$:
\begin{itemize}
\item[1.] Construct two arbitrary suitable immersions $S_1$ and $S_2$ with correct $J^+$-invariants.
\item[2.] Use the moves \textbf{M1} and \textbf{M3} and add the necessary exterior loops to the single immersions $S_1$ and $S_2$. Note that this action does not change their $J^+$-invariant. 
\item [3.] Push the exterior loops added \textbf{M3} into each other and the outer branch of $S_2$ into the exterior loops of $S_1$ added because of \textbf{M1}. This leads to the necessary amount of direct tangencies in the calculation of $J^{2+}(\overline{K})$ and $J^{2+}(K)$.
\end{itemize}
Case (ii): $z < x + y$ and $\overline{z} > 2x+2y - z$:
\begin{itemize}
\item[1.] Use move \textbf{M2} and construct $S_1$ including a suitable $L_i$ and an arbitrary immersion $S_2$ with correct $J^+$-invariants.
\item[2.] Push $S_2$ into the $L_i$-part of $S_1$ until the pair has the correct invariant $J^{2+}(K)$.
\item[3.] Use move \textbf{M3} to add the necessary amount of exterior loops to $S_1$ and $S_2$ and push them into each other until the pair has the correct invariant $J^{2+}(\overline{K})$.
\end{itemize}
Case (iii): $\overline{z} < x + y$ and $z > 2x+2y - \overline{z}$:
\begin{itemize}
\item[1.] Use move \textbf{M4} and construct an immersion $S_1$ including an exterior loop with suitable $L_i$ and an arbitrary immersion $S_2$ with correct $J^+$-invariants.
\item[2.] Push $S_2$ into the relevant part of $S_1$ until the pair has the correct invariant $J^{2+}(\overline{K})$. 
\item [3.] Use move \textbf{M1} to add the necessary amount of exterior loops to $S_1$ and $S_2$ and push them into each other until the pair has the correct invariant $J^{2+}(K)$.  
\end{itemize} \hrule
\ \\
\ \\
\begin{figure}[h]
\renewcommand*\figurename{Figure}
\begin{center}
\includegraphics[width=12cm]{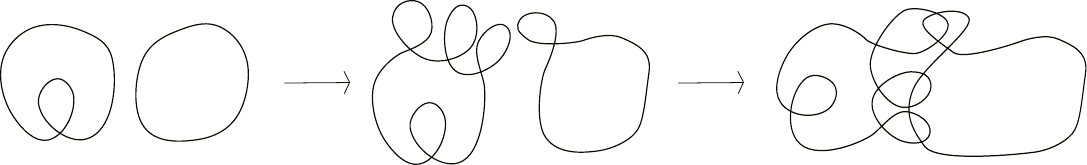}
\caption{Construction of a suitable pair of immersions with $J^+$-invariants $x=-2$, $y=0$ and $J^{2+}$-invariants $z = 2$ and $\overline{z}=0$}
\label{ex1}
\end{center}
\end{figure}
\ \\
\textbf{Example 1}: Find a suitable pair of immersions for the values $x=-2$, $y=0$, $z = 2$ and $\overline{z}=0$. Note that $z$, $\overline{z} \geq x + y = -2$, and proceed as given in \textit{Algorithm 1}.  
\begin{itemize}
\item [1.] Choose the immersions $K_2$ and $K_1$ as suitable immersions for $S_1$ and $S_2$. 
\item[2.] Calculate $d = z - x - y = 2 \cdot 2$ and $\overline{d} = 2 \cdot 1$ and therefore add according to \textbf{M1} two exterior loops to $S_1$ and according to \textbf{M3} one exterior loop to $K_2$ and $K_1$.
\item[3.] According to \textbf{M1} push the outer branch (not a loop) of $K_1$ into two exterior loops of $K_2$ and the left over exterior loops into each other. \end{itemize}
The pair went through two direct tangencies in the case of $J^{2+}(K)$ and one direct tangency in the case of $J^{2+}(\overline{K})$.
\begin{figure}[h]
\renewcommand*\figurename{Figure}
\begin{center}
\includegraphics[width=12cm]{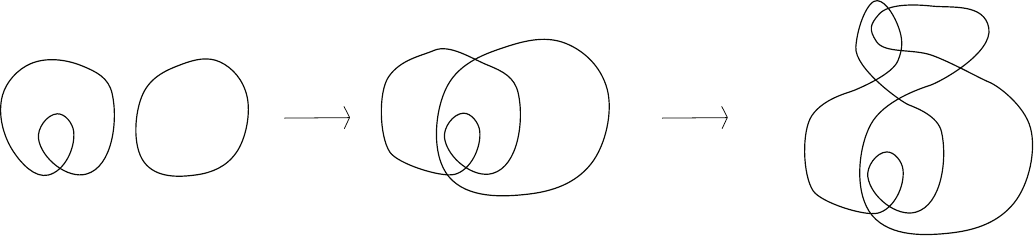}
\caption{Construction of a suitable pair of immersions with $J^+$-invariants $x=-2$, $y=0$ and $J^{2+}$-invariants $z = -4$ and $\overline{z}=4$}
\label{ex2}
\end{center}
\end{figure}\ \\
\textbf{Example 2}: Find a suitable pair of immersions for the values $x=-2$, $y=0$, $z =-4$ and $\overline{z}=4$. Note that $z \leq x + y = -2$, and $\overline{z} > 2x + 2y - z  = 0$, and proceed as given in \textit{Algorithm 1}.  
\begin{itemize}
\item [1.] Choose the immersions $K_2$ and $K_1$ as suitable immersions for $S_1$ and $S_2$.
\item[2.] Calculate $d = z - x - y = -2$ and push $K_1$ into $K_2$ to create one negative direct tangency for the case of $J^{2+}(K)$. Thus we created a pair with $J^{2+}(K')= -4$ and $J^{2+}(\overline{K'})=2$. 
\item[3.] According to \textbf{M3} add one exterior loop to both immersions which are then pushed into each other to create one positive direct tangency in the case of $J^{2+}(\overline{K})$. \end{itemize} 
\begin{figure}[h]
\renewcommand*\figurename{Figure}
\begin{center}
\includegraphics[width=12cm]{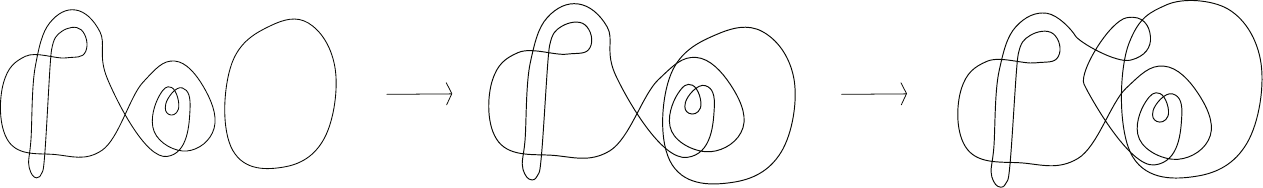}
\caption{Construction of a suitable pair of immersions with $J^+$-invariants $x=-4$, $y=0$ and $J^{2+}$-invariants $z = 4$ and $\overline{z}=-8$}
\label{ex3}
\end{center}
\end{figure}
\textbf{Example 3}: Find a suitable pair of immersions for the values $x=-4$, $y=0$, $z =4$ and $\overline{z}=-8$.\\
Note that $\overline{z} \leq x + y = -4$, and $z > 2x + 2y - \overline{z}  = 0$, and proceed as given in \textit{Algorithm 1}.
\begin{itemize}
\item [1.] Choose the immersions as visualized in the Figure \ref{ex3}. An immersion with $J^+$-invariant equal to 2 combined with an exterior loop containing $L_2$ as immersion $S_1$ and $K_1$ as suitable immersions for  $S_2$. 
\item[2.] Push $K_1$ trough the $L_2$-part of $S_1$. Thus we constructed a pair with $J^{2+}(\overline{K'}) = -8$ and $J^{2+}(K')=2$. 
\item[3.] According to \textbf{M1} add one exterior loop to $S_1$ and push the outer branch of $K_2$ into that exterior loop  to create one positive direct tangency in the case of $J^{2+}(K)$. \end{itemize}
\ \\
\textit{Algorithm 1} can not construct pairs of immersions with invariants $J^{2+}(K)$- and $J^{2+}(\overline{K})$ satisfying 
\begin{center} \begin{tabular}{ll}$z = x+y - d$ and $\overline{z} = x + y + d$ & or \\
$\overline{z} = x+y - d$ and $z = x + y + d$ & for $d \geq 0$, \end{tabular} \end{center}
which is equivalent to 
\begin{center} \begin{tabular}{ll}$z < x + y$ and $\overline{z} = x + y - z$ & or\\
$\overline{z} < x + y$ and $z = x + y - \overline{z}$. &\end{tabular}\end{center} 
This case will be realized by \textit{Algorithm 2}, which constructs pairs of immersions in which one immersion encycles the other immersion completely. Before introducing the algorithm, we define the immersion $A_k$ for $k \in \mathbb{N}$ as shown in Figure \ref{A_k}.
\ \\
The $J^+$-invariant of $A_k$ is given through:
\begin{center} $J^+(A_k) = 1 + 3 + 4 (k-1) - 2k - 0k = 4 + 4k - 4 - 2k = 2k$, \end{center} 
as with every $k$, the immersions $A_k$ has four more double points with index $0$ and two new components with winding number $1$ or $-1$.
\begin{lem} Let $K$ be a pair out the immersions $A_k$ for $k \in \mathbb{N}$ and $K_1$, whereas $K_1$ encircles $A_k$ completely without intersection. Then the $J^{2+}$-invariants are given through
\begin{center} $J^{2+}(K) = J^{2+}(\overline{K}) = 2k$. \end{center}
\begin{figure}[h]
\renewcommand*\figurename{Figure}
\begin{center}
\includegraphics[width=12cm]{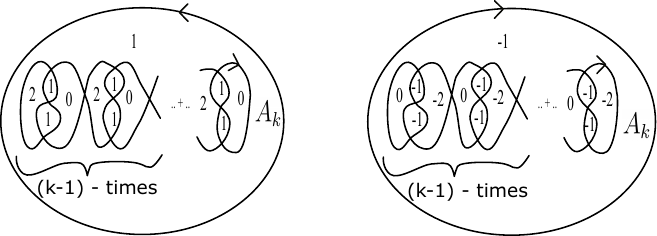}
\caption{Immersion $A_k$ completely encircled by $K_1$}
\label{A_k_en}
\end{center}
\end{figure}\end{lem}
\textit{Proof:} Figure \ref{A_k_en} shows $A_k$ encircled by $K_1$ with winding numbers for all components of $\mathbb{C}\setminus K$ for the calculation of $J^{2+}(K)$ and $J^{2+}(\overline{K})$. The $J^{2+}$-invariant is given through 
\begin{center}\begin{tabular}{lll}
$J^{2+}(K) = J^{2+}(\overline{K})$ & $=$ & $2 + 3 +4 (k-1) - (\pm 1)^2 - (\pm 2)^2 - 2(\pm1)^2 - (k-1) (2 (\pm 1)^2$ \\
 &  & $+ (\pm 2)^2) + 3 (\pm 1)^2+ (k-1)4(\pm 1)^2  +(\pm 1)^1 $\\
 & $=$ & $5 + 4k - 4 - 7 -6k +6 +3 + 4k -4 + 1$\\
 & $=$ & $2k$.\end{tabular}\end{center}$\hfill \Box$ \\
\ \\
\hrule
\ \\
\textbf{Algorithm 2:}\\
For given $x$, $y$, $z$ and $\overline{z}$ $\in 2\mathbb{Z}$ satisfying
\begin{center} (iv) $z = x + y - d$ and $\overline{z} = x + y + d$ for $d \in 2\mathbb{Z}$, \end{center} we construct a suitable pair with $J^+(S_1) = x$, $J^+(S_2) = y$, $J^{2+}(K)= z$ and $J^{2+}(\overline{K})=\overline{z}$. \\
For $d>0$, distinguish the following cases:\\
Case A: $x$ or $y \geq 2 - d$, choose without loss of generality $x \geq 2 - d$:
\begin{itemize}
\item[1.] Construct an immersion $S_1$ as connected sum of suitable $A_k$, $K_0$ and $K_{\frac{d}{2}}$ for $k=\frac{x + d - 2}{2}$, then $J^+(S_1) = x$. 
\item[2.] Choose an arbitrary immersion $S_2$ such that $J^+(S_2)=y$. 
\item[3.] The disjoint pair $K'$ has the invariants $J^{2+}(K')=J^{2+}(\overline{K'})=x + y$.
\item[4.] Now push $S_2$ through $S_1$ until $S_1$ is encircled once by $S_2$. Then the resulting pair has the invariants 
\begin{center} $J^{2+} (K)= x + y - d$ \\
$J^{2+}(\overline{K})= x + y + d$ \end{center} \end{itemize}
Case B: $x$, $y < 2 - d$, choose without loss of generality  $x < 2 - d$:
\begin{itemize}
\item[1.] Construct an immersion $S_1$ as $K_{1-\frac{x}{2}}$ with $k$ exterior loops, and $k = 1 - \frac{x+d}{2}$, then $J^+(S_1)= x$. 
\item[2.] Choose an arbitrary immersion $S_2$ such that $J^+(S_2)=y$. 
\item[3.] The disjoint pair $K'$ has the invariants $J^{2+}(K')=J^{2+}(\overline{K'})=x + y$.
\item[4.] Now push $S_2$ through $S_1$ until $S_1$ is encircled once by $S_2$. Then the resulting pair has the invariants 
\begin{center} $J^{2+}(K)= x + y - d$ \\
$J^{2+}(\overline{K})= x + y + d$ \end{center} \end{itemize} 
For $d<0$ (equivalent to $\overline{z}< x+y$) choose $S_2$ such that $J^+(S_2)=y$ and add one exterior loop. Let $K(S_1,S_2)$ be the pair of immersions with $S_1$ being completely encircled by the exterior loop of $S_2$. Then, $K$ has the suitable $J^{2+}$-invariants.\\ \hrule
\ \\
\ \\
Before going through the proof for \textit{Algorithm 2} there are two examples using the algorithm. \\
\begin{figure}[h]
\renewcommand*\figurename{Figure}
\begin{center}
\includegraphics[width=10cm]{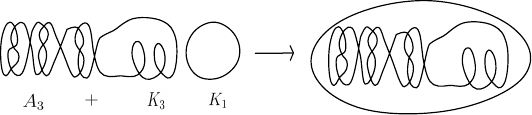}
\caption{Pair of immersions consisting of $K_1$ and an immersion as connected sum of $A_3$ and $K_3$. First pair is a disjoint pair, in the second pair $K_1$ encircles the other immersion completely.}
\label{Ex5}
\end{center}
\end{figure}\ \\
\textbf{Example 1:} We want to find a pair of immersions for $x=2$, $y=0$, $z = 8$ and $\overline{z}=-4$. Then the algorithm above can be used as 
\begin{center} \begin{tabular}{lll} $z$ & $= x + y + d = 2 + 0 +6 = 8$ & and\\
$\overline{z}$ & $= x + y - d = 2 + 0 - 6 = -4$ & with $d=6$ and $x \geq 2 - d = 2 - 6 = -4$ \end{tabular} \end{center}
\begin{itemize}
\item[1.] Construct an immersion $S_1$ as connected sum of $A_3$, $K_0$ and $K_{\frac{6}{2}}= K_3$, then \begin{center} $J^+(S_1) = J^+(A_2) + J^+(K_0)+ J^+(K_3) = 3\cdot 2 + 0 +(-4) = 2 = x$. \end{center}
\item[2.] Choose arbitrarily $K_1$ as $S_2$, \begin{center} $J^+(S_2) = J^+(K_1) = 0 = y$.\end{center}
\item[3.] The disjoint pair $K'$ fulfills $J^{2+}(K') = J^{2+}(\overline{K'})= x + y = 2$. 
\item[4.] Calculate the invariant of the pair in which $S_2$ completely encircles $S_1$, then
\begin{center} $J^{2+}(K) = 8 = z$ and $J^{2+}(\overline{K}) = -4 = \overline{z}$. \end{center}
\end{itemize}
\begin{figure}[h]
\renewcommand*\figurename{Figure}
\begin{center}
\includegraphics[width=10cm]{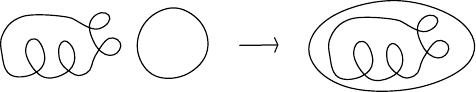}
\caption{Pair of immersions consisting of $K_3$ with two exterior loops and $K_1$. First pair is a disjoint pair, in the second pair $K_1$ encircles the other immersion completely.}
\label{Ex6}
\end{center}
\end{figure}\ \\
\textbf{Example 2:} We want to find a pair of immersions for $x=-4$, $y=0$, $z = -6$ and $\overline{z}=-2$. \textit{Algorithm 2} can be used as 
\begin{center} \begin{tabular}{lll}$z$ & $= x + y - d = -4+ 0 -2 = -6$ & and\\ 
$\overline{z}$ &  $= x + y - d = -4 + 0 + 2 = -2$ & with $d=2$ and $x < 2 - d = 2-2=0$ \end{tabular}\end{center}
\begin{itemize}
\item[1.] Construct an immersion $S_1$ as $K_{1-\frac{x}{2}} = K_{1-(-2)}=K_3$ and add $k = 1 - \frac{x+d}{2} = 1 - (-1) =2$ exterior loops, then \begin{center} $J^+(S_1) = J^+(K_3) = -4 = x$. \end{center}
\item[2.] Choose arbitrarily $K_1$ as $S_2$, \begin{center} $J^+(S_2) = J^+(K_1) = 0 = y$.\end{center}
\item[3.] The disjoint pair $K'$ fulfills $J^{2+}(K') = J^{2+}(\overline{K'})= x + y = -4$. 
\item[4.] Calculate the invariant of the pair in which $S_2$ completely encircles $S_2$, then 
\begin{center} $J^{2+}(K) = -4 - 2=-6 = z$ and $J^{2+}(\overline{K})=-4 + 2 = -2 = \overline{z}$.\end{center} \end{itemize}
\ \\
\textit{Proof for Algorithm 2:} Set $d>0$. For given $x$, $y$, $z$ and $\overline{z}$ which fulfill $z = x + y - d$ and $\overline{z} = x + y + d$ for a $d \in 2\mathbb{N}$, construct a suitable pair with $J^+{S_1} = y$, $J^+(S_2) = y$, $J^{2+}(K)= z$ and $J^{2+}(\overline{K})=\overline{z}$. \\
For $x$ with $x \geq 2 -d$ construct an immersion $S_1$ as connected sum of $K_{\frac{d}{2}}$, $K_0$ and $A_k$, whereas $k = \frac{x + d-2}{2}$. \\
Then the $J^+$-invariant is given through \begin{center} \begin{tabular}{lll}$J^+(S_1)$ & $=$ &  $J^+(K_{\frac{d}{2}}) +J^+(K_0) + J^+(A_k) $\\
 & $=$ & $2 - 2(\frac{d}{2})+0+ 2k$ \\
 & $=$ & $2 - 2(\frac{d}{2})+0+ 2(\frac{x+d-2}{2})$ \\
 & $=$ &  $ 2 - d + x + d - 2 $\\
 & $=$ &  $ x$\end{tabular}\end{center} 
Choose an arbitrary suitable immersion $S_2$ with $J^+(S_2) = y$ and then form a pair $K$ such that $S_1$ is once encircled by $S_2$ without intersection. \\
Then calculate the $J^{2+}$-invariants and use $x = 2k + 2 - d$:\begin{center}\begin{tabular}{lll}
$J^{2+}(K)$ & $=$ & $1 + 1 + n_{S_2}+ \frac{d}{2} -1 + 4k $\\
 & & $ - \sum_{C \subset \mathbb{C} \setminus S_2} (\omega_C(K))^2 - (\pm 2)^2 - (\frac{d}{2}-1) (\pm 3)^2 - k (2(\pm 1)^2 + (\pm 2)^2) $\\
 & & $ + \sum_{p \in S_2} (ind_p(K))^2 + (\frac{d}{2}-1)(\pm 2)^2 + k(4(\pm 1)^2)$)\\
 & & $ + (\pm 1)^2$\\
 & $=$ & $y + 2k + 2 - 2d$\\
& $=$ & $x + y -d$\end{tabular}\end{center}
For visualization check the example in Figure \ref{Ex5}.\\
\begin{center}\begin{tabular}{lll}
$J^{2+}(\overline{K})$ & $=$ & $1 + 1 + n_{S_2}+ \frac{d}{2} -1 + 4k $\\
 & & $- \sum_{C \subset \mathbb{C} \setminus S_2} (\omega_C(\overline{K}))^2 - 0 - (\frac{d}{2}-1)(\pm 1)^2 - k (2(\pm 1)^2 + (\pm 2)^2) $\\
 & & $ + \sum_{p \in S_2} (ind_p(\overline{K}))^2 + (\frac{d}{2}-1)( 0)^2 + k(4(\pm 1)^2) $\\
 & & $ + (\pm 1)^2$\\
 & $=$ & $y + 2 + 2k $\\
& $=$ & $x + y +d$\end{tabular}\end{center}
\ \\
For $x$ with $x <  2 -d$ construct an immersion $S_1$ as connected sum of $K_{1-\frac{x}{2}}$ and add $k$ exterior loops, whereas $k = 1- \frac{x}{2}-\frac{d}{2}$.\\
Then the $J^+$-invariant is given through \begin{center} \begin{tabular}{lll}$J^+(S_1)$ & $=$ & $J^+(K_{1-\frac{x}{2}})$\\
 & $=$ & $2 - 2 (1 - \frac{x}{2})  = x$\end{tabular}\end{center} 
Choose an arbitrary suitable immersion $S_2$ with $J^+(S_2) = y$ and then form a pair $K$ such that $S_1$ is completely encircled by $S_2$ without intersection. \\
Then calculate the $J^{2+}$-invariants:\begin{center}\begin{tabular}{lll}
$J^{2+}(K)$ & $=$ & $1 + 1 +  n_{S_2}+ (- \frac{x}{2}) + k $\\
 & & $ - \sum_{C \subset \mathbb{C} \setminus S_2} (\omega_C(K))^2 - (\pm 2)^2 - (-\frac{x}{2})(\pm 3)^2 - k(0)^2 $ \\
 & & $ + \sum_{p \in S_2} (ind_p(K))^2 + (-\frac{x}{2})(\pm 2)^2 + k (\pm 1)^2 $\\
 & & $ + (\pm 1)^2$\\
 & $=$ & $y - 2 + 2x + 2k$\\
 & $=$ & $y - 2 + 2x + 2(1 - \frac{x}{2}-\frac{d}{2})$\\
& $=$ & $x + y -d$\end{tabular}\end{center}
\begin{center}\begin{tabular}{lll}
$J^{2+}(\overline{K})$ & $=$ & $1 + 1 +  n_{S_2}+ (- \frac{x}{2}) + k $\\
 & & $ - \sum_{C \subset \mathbb{C} \setminus S_2} (\omega_C(\overline{K}))^2 - (-\frac{x}{2})(\pm 1)^2 - 4k$ \\
 & & $ + \sum_{p \in S_2} (ind_p(\overline{K}))^2 + (-\frac{x}{2})(0)^2 + k$ \\
 & & $ + (\pm 1)^2$\\
 & $=$ & $y + 2 - 2k$\\
 & $=$ & $y + 2 - 2(1 - \frac{x}{2}-\frac{d}{2})$\\
& $=$ & $x + y + d$\end{tabular}\end{center}
For visualization check the example in Figure \ref{Ex6}. The proof for the case $d<0$ works analogously.$\hfill \Box$ \\
\ \\
\textit{Proof of Theorem D:} The introduced \textit{Algorithms 1} and \textit{2} construct suitable pairs of immersions for the cases\begin{center}
\begin{tabular}{lll}
\textit{Algorithm 1} & (i) & $z, \overline{z} \geq x + y$\\
 & (ii) & $z < x + y$ for and $\overline{z} > 2x+2y - z$ or\\
 & (iii)& $\overline{z} < x + y$ for and $z > 2x+2y - \overline{z}$ \\
\textit{Algorithm 2} & (iv) & $z = x + y - d$ and $\overline{z} = x + y + d$ for a $d \in 2\mathbb{Z}$ \end{tabular} \end{center}
It remains to show that those cases cover all possibilities of combinations of $x$, $y$, $z$ and $\overline{z}$ with $z + \overline{z} \geq2 x + 2y$.\\
Let $x$, $y$, $z$ and $\overline{z}$ $\in 2 \mathbb{Z}$ be arbitrary satisfying the above equation but not satisfying the above cases (i), (ii), (iii) and (iv). Because of (i) either 
\begin{center} $z$ or $\overline{z} < x + y$,\end{center}
without loss of generality choose $z < x + y$. As (ii), (iii) and (iv) are also not satisfied $\overline{z}$ satisfies \begin{center}$\overline{z} < 2x + 2y - z$.\end{center}
This leads to
\begin{center} $z + \overline{z} < z + 2x+ 2y - z = 2x + 2y$\end{center}
which is a contradiction to the initial conditions. $\hfill \Box$\\
\begin{figure}[h]
\renewcommand*\figurename{Figure}
\begin{center}
\includegraphics[width=10cm]{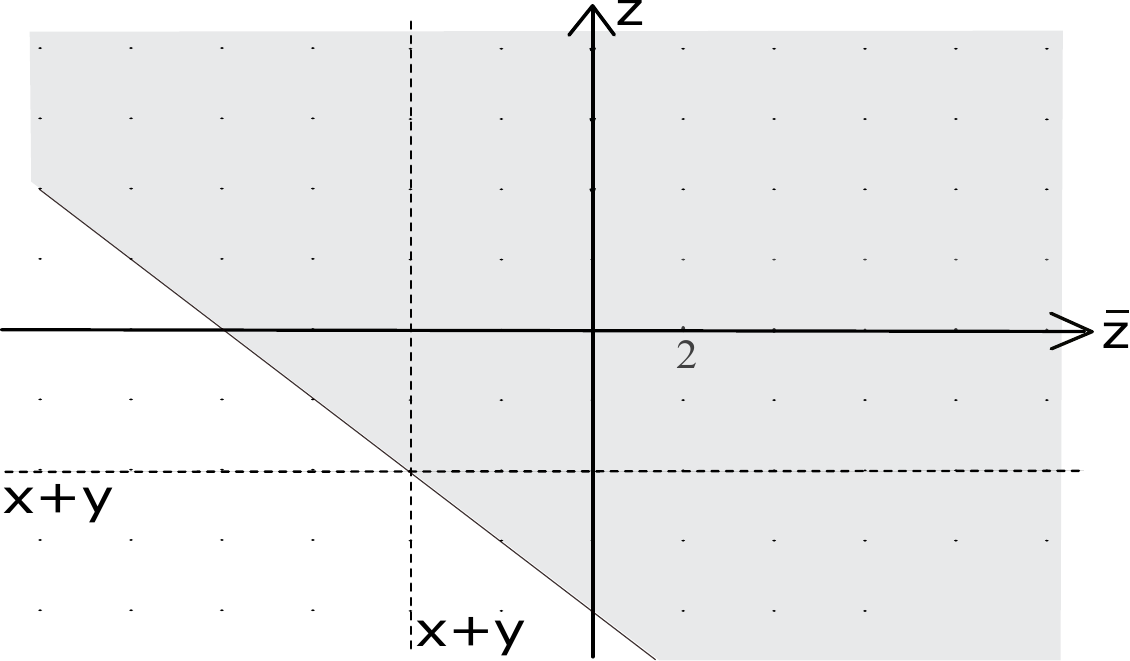}
\caption{Possible combinations of $J^{2+}$-invariants $z$ and $\overline{z}$ depending on the sum of $J^+$-invariants $x + y$}
\label{solutions}
\end{center}
\end{figure}\\
Figure \ref{solutions} shows which combinations of $z$ and $\overline{z}$ exist for a given sum $x + y$. In the graph the possible combinations of invariants $z$ and $\overline{z}$ are placed on or above the linear graph going through ($x+y$, $x+y$). The grey part can be realized with \textit{Algorithm 1}, whereas \textit{Algorithm 2} constructs pairs with invariants visualized through the linear graph going through ($x+y$, $x+y$), the exact point ($x+y$, $x+y$) can be realized as a disjoint pair of immersions. 
 \subsection{Interpretation as $J^{2-}$-invariant}
In [1], Arnold did not only invent the $J^+$-invariant but also the $J^-$-invariant which is sensible to inverse tangencies under homotopy and the $St$ Strangeness which is sensible to triple points. The $J^-$-invariant is unique up to an additive constant for generic immersions of fixed index whose value remains unchanged while the immersed curve experiences a direct tangency or a triple point crossing, but decreases by $2$ under a positive crossing of an inverse tangency. It is defined for the standard immersions $K_j$ through
\begin{center} $J^-(K_0)=-1$ and $J^-(K_j)=-3(j-1)$ for $j \geq 1$. \end{center}
Arnold shows the following relation $J^+ - J^- = n$, therefore, we calculate $J^{-}$ using Viro's formula.
\begin{align*}
      J^-(S) &  = J^+(S) - n = 1 + n - \sum_{C\subset \mathbb{C}\setminus S} \omega_C(S)^2 + \sum_p ind_p(S)^2 -n \\
    & = 1  - \sum_{C\subset \mathbb{C}\setminus S} \omega_C(S)^2 + \sum_p ind_p(S)^2\\  
\end{align*}
\begin{theo}
    The $J^-$-invariant decreases (increases) by the value of $2$ when going positively (negatively) through inverse tangency and is constant under direct tangency and triple points.
\end{theo}
\textit{Proof:} First, we calculate the $J^-$-invariants for the standard immersions: \\
Set $j=0$: $J^-(K_0)=1 - 1^2 - (-1)^2 + 0^2 = -1$\\
Set $j\geq 1$: $J^-(K_j) = 1 - 1^2 - (j-1) \cdot 2^2 + (j -1) \cdot 1^2 = - 4(j-1)+(j-1) = -3 (j-1)$\\
Now, we check the behavior under the critical scenarios analogously to the poof of Viro's formula to calculate $J^+$-invariant. When going through a direct or inverse tangency, the change of the invariant is given through 
\begin{align*}
    \Delta J^-(S)   & = - \Delta \sum_{C\subset \mathbb{C}\setminus S} \omega_C(S)^2 + \Delta \sum_p ind_p(S)^2\\
                    & = z^2 + s^2 + 2 (\frac{1}{4}(x+y+z+s))^2,
\end{align*}
as visualized in Figure \ref{Viro_Bew4}.
Considering the correct winding numbers of the components in the direct setting, we receive
\begin{equation*}
    \Delta J^-(S) = -2z^2 + 2z^2 = 0
\end{equation*}
and the winding numbers in the inverse setting, we get
\begin{equation*}
    \Delta J^-(S) = -z^2 - (z+2)^2 + 2(z+1)^2 = - 2
\end{equation*}
When going through a triple point, we receive 
\begin{align*}
\triangle {J}^{-}(S) &  = - \triangle \sum \omega_C (S)^2 + \triangle \sum ind_p(S)^2 \nonumber \\
 & = - \left( \omega_{T_{\mathrm{new}}}(S)^2 - \omega_{T_{\mathrm{old}}}(S)^2 \right) + \triangle ind_A(S)^2 + \triangle ind_B(S)^2 + \triangle ind_C(S)^2
 & = 0.
\end{align*}
where $T$ describes the old and new triangle formed by the double points $A$, $B$ and $C$ of $S$, as visualized in Figure \ref{Viro_Bew6}. $\hfill\Box$\\
\ \\
In view of the formulation of $J^{2+}$, we define
\begin{equation*}
    J^{2-}(K):= 2 - \sum_{C\subset \mathbb{C}\setminus S} \omega_C (S)^2 + \sum_p ind_p(S)^2 + u(K)^2.
\end{equation*}
\begin{cor}
    The $J^{2-}$-invariant of a disjoint pair $K(S_1,S_2)$ is given through $J^-(S_1) + J^-(S_2)$.
\end{cor}
\textit{Proof:} This follows immediately due to the construction of $J^{2-}$. $\hfill \Box$\\
\begin{theo}
    The $J^{2-}$-invariant for pairs of immersions decreases (increases) by the value of $2$ when going through positive (negative) inverse tangency and does not change through direct tangency and triple points. 
\end{theo}
\textit{Proof:} The proof works analogously to the proof of the $J^{2+}$-invariant, in which we check the behavior under the six critical scenarios. $\hfill\Box$\\
\begin{prop}
 Let $K(S_1,S_2)$ be an oriented pair of immersions, the relation between $J^{2+}$ and $J^{2-}$-invariants is given through 
\begin{equation*}
    J^{2-}(K) = - J^{2+}(\overline{K}) + J^+(S_1) + J^+(S_2) + J^-(S_1) + J^-(S_2)
\end{equation*}   
\end{prop}
\textit{Proof:} We can easily see that the formulas correspond for disjoint pairs of immersions. It remains to show that the formula behave equally under the critical scenarios,
\begin{equation*}
    \Delta J^{2-}(K) = - \Delta J^{2+}(\overline{K}) + \Delta J^+(S_1) + \Delta  J^+(S_2) + \Delta  J^-(S_1) + \Delta  J^-(S_2)
\end{equation*}
We first go through all critical scenarios during single-homotopy. Let $S_1$ be the immersion going through a critical scenario
\begin{center}\begin{tabular}{lll}
Direct tangency & $\Delta J^{2-}(K)$ & $= - \pm 2 + (\pm 2) + 0+0+0 = 0$\\
Inverse tangency & $\Delta J^{2-}(K)$ & $= - 0 + 0+0+(\mp 2) + 0 = \mp 2$\\
Triple points  & $\Delta J^{2-}(K)$ & $= - 0 + 0 + 0+0+0 = 0$\\
\end{tabular}\end{center}
When going through critical scenarios during an unravel-homotopy, we get
\begin{center}\begin{tabular}{lll}
Direct tangency & $\Delta J^{2-}(K)$ & $= - 0 + 0 + 0+0+0 = 0$\\
Inverse tangency & $\Delta J^{2-}(K)$ & $= - (\pm 2) + 0+0+0 + 0 = \mp 2$\\
Triple points  & $\Delta J^{2-}(K)$ & $= - 0 + 0 + 0+0+0 = 0$\\
\end{tabular}\end{center}
$\hfill\Box$\\

\section{Generalization for $n$ immersions}
Arnold's $J^+$-invariant can be extended to a triple of oriented immersions $(S_1,S_2,S_3)$ defined as 
\begin{equation*}
(S_1,S_2,S_3) : S^1 \sqcup S^1 \sqcup S^1 \rightarrow \mathbb{C}
\end{equation*}
by the following formula
\begin{equation*}
    J^{3+}(S_1,S_2,S_3) = J^{2+}(S_1,S_2) + J^{2+}(S_2,S_3) + J^{2+}(S_1,S_3)  - J^+(S_1) - J^+(S_2) - J^+(S_3) 
\end{equation*}
and extended $\forall n \geq 3$ for links of $n$ oriented immersions $(S_1,S_2,..,S_n)$ by
\begin{equation*}
    J^{n+} (S_1, S_2,..,S_n) = \sum_{i<j} \Bigl( J^{2+}(S_i, S_j)-J^+(S_i) - J^+(S_j)\Bigr) + \sum_{i=1}^n J^+(S_i) .
\end{equation*} 
Analogously to the definition of $J^{n+}$, we define $\forall n \geq 3$
\begin{equation*}
   J^{n-}(S_1, S_2, .. ,S_3) := \sum_{i<j}\Bigl( J^{2-} (S_i,S_j)-J^-(S_i)-J^-(S_j)\Bigr) + \sum_{i=1}^n J^-(S_i).
\end{equation*}
\begin{theo}[Theorem E]
    The $J^{n+(-)}$-invariant for links of $n$ oriented immersions $S_1,S_2,..,S_n$ changes by the value of $2$ when going through direct (inverse) tangency and remains constant under inverse (direct) tangency and triple points. 
\end{theo}
\textit{Proof:} During a homotopy of a system of $n$ immersions seven different critical scenarios can occur:\\
\begin{center}\begin{tabular}{ll}
1. & direct tangency during single-homotopy\\
2.\& 3.& inverse tangency and triple points during single-homotopy\\
4. & direct tangency during unravel-homotopy\\
5.\& 6. & inverse tangency and triple points during unravel-homotopy including two immersions\\
7. & triple point during unravel-homotopy including three immersions
\end{tabular}\end{center}
Let us first consider the $J^{n+}$-invariant. As the $J^+$- and $J^{2+}$-invariants do not change under triple points and inverse tangencies, $J^{n+}$ also does not change under these scenarios. It remains to show that $J^{2+}$ changes by the value of $2$ when going through direct tangencies, namely the scenarios $1.$ and $4.$ Let $S_1$ be the transformed immersion in $1.$. Then only the $J^{2+}$-invariants for pairs $K(S_i,S_j)$ for $i=1$ and $\forall j \in [2,n]$ and the $J^+$-invariant of $S_1$ do change by the value of $\pm2$. Hence, we get\\ 
\begin{center}
 \begin{tabular}{ll}
  $\Delta J^{n+}$ & $ = \Delta \sum_{i<j} \Bigl( J^{2+}(S_i, S_j)-J^+(S_i) - J^+(S_j)\Bigr) + \Delta \sum_{i=1}^n J^+(S_i)$\\
& $ = \sum_{j=2}^n \Bigl( \Delta J^{2+}(S_1, S_j)-\Delta J^+(S_1) - \Delta J^+(S_j)\Bigr) + \Delta J^+(S_1)$\\
    & $ =  \sum_{j=2}^n \Bigl( \pm 2 - (\pm 2) - 0 \Bigr) + (\pm2) = \pm 2$ \end{tabular}\end{center}
Let $S_1$ and $S_2$ be the immersions going through a direct tangency in scenario $4.$ Then the $J^{2+}$-invariant of the pair $K(S_1,S_2)$ changes by the value $2$, whereas all other $J^{2+}$-invariants and all $J^+$-invariants remain constant. \begin{center} \begin{tabular}{ll}
  $\Delta J^{n+}$ & $ = \Delta \sum_{i<j} \Bigl( J^{2+}(S_i, S_j)-J^+(S_i) - J^+(S_j)\Bigr) + \Delta \sum_{i=1}^n J^+(S_i)$\\
  & $ = \Bigl( \Delta J^{2+}(S_1, S_2)-\Delta J^+(S_1) - \Delta J^+(S_2)\Bigr) + \Delta \sum_{i=1}^n J^+(S_i)$\\
  & $ =  (\pm 2 - 0 - 0 \Bigr)  + 0 = \pm 2$ \end{tabular}\end{center} 
The proof for $J^{n-}$ works analogously.$\hfill\Box$\\
\begin{defi} A link of $n$ immersions $(S_1, .. , S_n)$ is called \textbf{disjoint}, if $S_i \subset C^* \in \mathbb{C}\setminus \{S_j\}$ $\forall i,j \in [1,n]$, $ i \neq j$ and $C^*$ being the unbounded component in $\mathbb{C}\setminus \{S_j\}$.\end{defi}
\begin{cor} Let $(S_1, .. , S_n)$ be a disjoint link of $n$ immersions, then $J^{n\pm}= \sum_{i=1}^n J^{\pm}(S_i)$. \end{cor}
\textit{Proof:} This follows immediately as for a disjoint pair $J^{2\pm}(K) = J^\pm(S_1) + J^\pm(S_2)$. $\hfill \Box$\\

\newpage
\section*{References}
\begin{small}
 $[1]$ V.I. Arnold, \textit{Plane curves, their invariants, perestroikas and classifications}, Adv. Soviet Math. \textbf{21} (1994).  \\
$[2]$ K. Cieliebak, U. Frauenfelder and O. van Koert, \textit{Periodic orbits in the restricted three-body problem and Arnold’s $J^+$-invariant}, Regul. Chaotic Dyn. \textbf{22} (2017) 408–434.\\
$[3]$ U. Frauenfelder and O. van Koert,   \textit{The restricted three-body problem and holomorphic curves}, Pathways in Mathematics, Birkhäuser Basel (2018).\\
$[4]$ H. Haeussler, \textit{The $J^{2+}$-invariant for pairs of generic immersions}, arXiv: 2104.04349 (2021).\\ 
$[5]$ J. Kim and S. Kim,  \textit{$J^+$-like invariants of periodic orbits of the second kind in the restricted three body problem}, Journal of Topolgz and  Analysis \textbf{12} (2017).\\
$[6]$ O. Viro,  \textit{Generic immersions of circle to surfaces and complex topology of real algebraic curves}, Amer. Math. Soc. Transl. Ser. 2. \textbf{173} (1996) 231-252 \\
$[7]$ H. Whitney, \textit{On regular closed curves in the plane}, Compositoo Mathematica 4 (1937) 276-284. \\   
\end{small}

\end{document}